\documentclass[12pt]{article}
\usepackage{epsfig}
\usepackage{graphics}
\usepackage{amsmath,amsthm,amssymb,amsbsy,graphicx,epsfig,citesort}
\numberwithin{equation}{section}
\providecommand{\C}[1]{\mathcal{#1}}
\providecommand{\D}[1]{\mathbb{#1}}

\newcommand{\dist}{\mathop{\rm dist}}
\renewcommand{\div}{\mathop{\rm div}}

\theoremstyle{plain}
\newtheorem{thm}{Theorem}%[section]

\newtheorem{lem}[thm]{Lemma}
\newtheorem{prop}[thm]{Proposition}
\theoremstyle{definition}
\newtheorem{Con}[thm]{Conjecture}

\newtheorem{rem}[thm]{Remark}
\newtheorem*{assumps*}{Assumptions}
\newtheorem*{nots*}{¡®§­ ç¥­¨ï}
\newtheorem*{rem*}{ ¬¥ç ­¨¥}
\def \ve{\varepsilon}
\begin{document}
%\inputencoding{cp866}
%\begin{large}
%\setlength{\labelwidth}{40pt}

\begin{center}
\begin{Large}
Solutions with Vortices of a Semi-Stiff Boundary Value Problem
for the Ginzburg-Landau Equation
\end{Large}
\end{center}
\bigskip
\begin{center}
{\bf  L. Berlyand}\footnote{Partially supported by
NSF grant DMS-0708324}
\end{center}

\noindent
Department of Mathematics, The Pennsylvania
State University, University Park, PA 16802, USA,

\noindent
email: berlyand@math.psu.edu

\begin{center}
and
\end{center}

\begin{center}
{\bf V. Rybalko}\footnote{Partially supported by
the Akhiezer fund and by Shapiro fund at Penn State University}
\end{center}

\noindent
Mathematical Division, Institute for Low Temperature Physics and
Engineering, 47 Lenin Ave., 61103
Kharkov, Ukraine

\noindent
email:vrybalko@ilt.kharkov.ua

\bigskip
\begin{abstract}
We study solutions of the 2D Ginzburg-Landau equation
$$
-\Delta u+\frac{1}{\ve^2}u(|u|^2-1)=0
$$
subject to "semi-stiff"  boundary conditions: the Dirichlet
condition for the modulus, $|u|=1$, and the homogeneous Neumann
condition for the phase. The principal result of this work shows
there are stable solutions of this problem with zeros (vortices),
which are located near the boundary and have bounded energy in the
limit of small $\ve$.    For the Dirichlet bondary condition
("stiff" problem), the existence of  stable solutions with vortices,
whose energy blows up as $\ve \to 0$, is well known. By contrast,
stable solutions with vortices are not established in the case of
the homogeneous Neumann ("soft") boundary condition.
%(nonexistence is proved for simply connected domains).

In this work, we develop a variational method which allows one to construct
local minimizers of the corresponding Ginzburg-Landau
energy functional. We introduce an approximate bulk degree as the key
ingredient of this method, and, unlike the standard degree over the curve,
it is preserved in the weak $H^1$-limit.
\end{abstract}

%by studying a constrained minimization problem
%with the following two constraints. The first constraint consists in
%prescribing topological degrees on the boundary. The second one is
%an interval constraint on

%The proof is based on the introduction of a notion of the
%bulk winding number which is stable with respect to the weak $H^1$
%convergence unlike the standard degree boundary conditions.

\section{Introduction and Main Results}
In this work, we study solutions of the Ginzburg-Landau (GL) equation
\begin{equation}
-\Delta u+\frac{1}{\ve^2}u(|u|^2-1)=0 \quad\text{in}\quad A,
\label{equation}
\end{equation}
%where $A$ is a bounded multiply connected domain in $\D{R}^2$,
where $\ve$ is a positive parameter (the inverse of the
GL parameter $\kappa=1/\ve$), $u$ is a complex-valued
($\D{R}^2$-valued) map, and $A$ is a smooth, bounded, multiply connected domain in $\D{R}^2$.
For
simplicity, hereafter we assume $A$ is an annular type (doubly
connected) domain of the form
$A=\Omega\setminus\overline{\omega}$, where $\Omega$ and $\omega$ are
simply connected smooth domains and $\overline{\omega}\subset\Omega\subset\D{R}^2$.

Equation (\ref{equation}) is the Euler-Lagrange PDE
%partial differential equation
corresponding to the energy functional
\begin{equation}
E_\ve(u)=\frac{1}{2}\int_A|\nabla u|^2{\rm d}x
+\frac{1}{4\ve^2}\int_A (|u|^2-1)^2{\rm d}x.
\label{P1}
\end{equation}
Equations of this type
%(\ref{equation})
arise, e.g., in models of superconductivity and superfluidity.  Additionally,
%various phase transition problems occurring e.g. in superfluidity or superconductivity.
(\ref{equation}) is viewed as a complex-valued version of the
Allen-Cahn model for phase transitions \cite{Sreview}.

Solutions of (\ref{equation}) subject to the Dirichlet boundary
condition, $u=g$ on $\partial A$ with fixed $S^1$-valued boundary
data $g$, have been extensively studied in the past decade. Special
attention has been paid to solutions with isolated zeros (vortices).
In contrast with the Dirichlet problem, in the case of the
% (\ref{equation}) with
homogeneous Neumann boundary condition, solutions are typically
vortexless; in particular, stable solutions with vortices have not
been established.

%received much attention in the literature as well.
%It corresponds to finding critical points of the energy (\ref{P1}) in the entire space
%$H^1(A;\D{R}^2)$.

This work is devoted to solutions
of (\ref{equation}) subject to the ``semi-stiff"
boundary conditions
\begin{equation}
|u|=1\ \text{and}\ u \times \frac{\partial u}{\partial \nu}=0 \quad\text{on}\ \partial A.
\label{boundary}
\end{equation}
These boundary conditions are intermediate between Dirichlet and
Neumann in the following sense: any solution $u\in H^1(A;\D{R}^2)$
of (\ref{equation}, \ref{boundary}) is sufficiently regular
\cite{BM1}, so it can be written as $u=|u|{\rm e}^{i\psi}$ (locally)
near the boundary. Then (\ref{boundary}) means the Dirichlet
boundary condition is prescribed for the modulus, $|u|=1$ on
$\partial A$, and the Neumann condition is prescribed for the phase,
$\frac{\partial \phi}{\partial \nu}=0$ on $\partial A$.

Problem (\ref{equation}, \ref{boundary}) is equivalent to finding
critical points of the energy functional (\ref{P1}) in the space
\begin{equation}
\label{definitionofclassJ}
\C{J}=\{u\in H^1(A;\D{R}^2);\ |u|=1\ \text{a.e. on}\ \partial A\}.
\end{equation}
Our main objective is to study the existence of {\it stable}
solutions of (\ref{equation}, \ref{boundary}) with  vortices. Since
the problem is time independent, stable solutions are defined as
(local) minimizers of (\ref{P1}) in $\C{J}$. In other words, we are
interested in whether the model (\ref{equation}, \ref{boundary})
%, which is intermediate between
%the classical Dirichlet and Neumann ones,
stabilizes vortices similarly to Dirichlet problem or does not
stabilize vortices analogously to Neumann problem.  The boundary conditions (\ref{boundary})
are not well studied, and this work, along with studies
\cite{BV,GB,BM1,BM2,BGR} reveals their distinct features, described later in the introduction.

%While
%vortexless solutions of (\ref{equation}),(\ref{boundary}) have been obtained in \cite{GB} %and \cite{BM1},\cite{BM2}, the existence of solutions with nontrivial vortex structure
%is established in this work for the first time.

%which have bounded energy in the limit of small $\ve$.

%If a solution $u$ of (\ref{equation}) has an isolated zero, around which
%the {\it winding number},
%or {\it degree} of $u/|u|$ being nonzero, such a point is called {\it vortex}.
%The existence and stability of solutions with vortices crucially depends on the boundary %conditions imposed on the boundary.
Let us briefly review the existing results for the Dirichlet and
Neumann boundary value problems for equation (\ref{equation}). The
first results on the existence of stable solutions with vortices for
Dirichlet problem were obtained in \cite{EMT-Q,FP}. Stable solutions
of (\ref{equation}) with vortices were obtained and studied in
\cite{BBH} for star-shaped domains and prescribed $S^1$-valued
boundary data with nonzero topological degree. In \cite{BBH}, the
limiting locations (as $\ve\to 0$) of vortices of globally
minimizing and other solutions (if they exist) are described by
means of a {\it renormalized energy}. Subsequently, these results
were generalized for multiply connected domains in \cite{Str}. The
existence of locally minimizing and minmax solutions was established
first in \cite{L} and \cite{LL}, then in more detail and generality
in \cite{PR} (see also \cite{AB,DKM}). We refer the reader to
\cite{B?}, and references therein, for the Dirichlet problem's
various results. As previously mentioned, only vortexless stable
solutions of (\ref{equation}) with the homogeneous Neumann boundary
condition in 2D are known. Moreover, all locally minimizing
solutions are constant maps if $A$ is convex \cite{JM}, or simply
connected and $\ve$ is small \cite{S1}. The existence of nonconstant
(but vortexless) locally minimizing solutions is established in
\cite{JMZ} and \cite{A}. In the recent work \cite{DKM}, a general
result for the existence of (nonminimizing) solutions with vortices
was obtained. Similarly to the Dirichlet problem, these solutions
with vortices have energy that blows up as $\ve\to 0$.

%In summary one can
%say that the model (\ref{equation}) with the Dirichlet boundary condition
%(loosely speaking) stabilizes
%vortices while one with the Neumann boundary condition allows no stable solutions
%with vortices (at least such solutions are not known yet).
Equation (\ref{equation}) (functional (\ref{P1})) is usually referred to as a
simplified GL
model (without magnetic field). There is a large body of
mathematical literature on the general GL model with a magnetic field
(e.g., \cite{SS1,BP,ABr,SR,JS}).
Since (\ref{equation}) is obtained from the general GL energy by
setting the magnetic field to zero, it describes persistent
currents in a 2D cross-section of a cylindrical superconductor
(or in a 2D film). It was observed in \cite{BBH}, the degree
of the boundary data on connected components of $\partial A$
creates the same type of "quantized vortices" as  a magnetic field
in type II superconductors or as angular rotation in superfluids.
Despite a relatively simple form of equation (\ref{equation}), it
leads to a
deep analysis of properties of its solutions similar to other
fundamental PDE's
% partial differential equations
in mathematical physics.

The boundary conditions (\ref{boundary})
model, e.g., the surface of a superconductor coated with a
high temperature superconducting thin film \cite{Ar}.  Generating a
mathematical model of persistent currents in such a superconductor, then
amounts to finding critical points of functional (\ref{P1}) in the space
$\C{J}$, when $u=|u|e^{i\psi}$ on the boundary and $|u|=1$, while
the phase $\psi$ is "free".

Boundary conditions (\ref{boundary}) appeared in recent studies
\cite{BV,GB,BM2} of the minimization problem for the GL functional
(\ref{P1}) among maps from $\C{J}$ with prescribed degrees on the
connected components of the boundary. The minimization of the energy
(\ref{P1}) in $\C{J}$ produces only constant solutions of
(\ref{equation}, \ref{boundary}), similar to the case of the Neumann
problem, which corresponds to finding critical points of (\ref{P1})
in the entire space $H^1(A;\D{R}^2)$. An obvious way of producing
critical points with vortices  is to impose two different degrees
$q\not=p$ on $\partial\Omega$ and $\partial\omega$. That is, to
consider the minimization of $E_\ve(u)$ in the set
$\C{J}_{pq}\subset\C{J}$, where
\begin{equation}
\label{definitionofclassJpq}
\C{J}_{pq}:=\{u\in \C{J};\ {\rm deg}(u,\partial\omega)=p,{\rm deg}(u,\partial\Omega)=q. \}
\end{equation}
Recall that the degree (winding number) of a map $u\in H^{1/2}(\gamma,S^1)$ on
$\gamma$ (where $\gamma$ is either $\partial \omega$ or $\partial
\Omega$) is an integer given by the classical formula (cf., e.g.,
\cite{B})
\begin{equation}
\label{def_degree}
\deg(u,\gamma)=\frac{1}{2\pi}\int_\gamma u\times\frac{\partial u}{\partial \tau}{\rm d}\tau,
\end{equation}
where the integral is understood via $H^{1/2}-H^{-1/2}$ duality, and
$\frac{\partial }{\partial \tau}$ is the tangential derivative with
respect to the counterclockwise orientation of $\gamma$. (Throughout
the paper we assume the same orientation of $\partial \omega$ and
$\partial \Omega$.) Note that $\C{J}_{pq}$ are connected components
of $\C{J}$ (see \cite{B}).

Simple topological considerations imply that critical points from
$\C{J}_{pq}$ must have at least $|p-q|$ (with multiplicity)
vortices. We emphasize that the existence of such critical points is
far from obvious. For example (see Section \ref{section2}),
there are no global minimizers of $E_\ve(u)$ in $\C{J}_{01}$ and
the weak limits of minimizing sequences do not belong to
$\C{J}_{01}$. This simple example illustrates an important
property of the sets $\C{J}_{pq}$, which is crucial for our
consideration: {\it these sets are not weakly $H^1$-closed}, since
the degree at the boundary may change in the limit. On the other
hand, the results of this work show that when $p=q$ and there is no
topological reason for vortices to appear, local minimizers
typically do have vortices.

%\subsection{Main results}
\smallskip

As mentioned above, the vortex structure of solutions of
(\ref{equation}) with Dirichlet and Neumann boundary conditions is
well studied. In contrast, only vortexless solutions of the
semi-stiff problem (\ref{equation}), (\ref{boundary}) were found
\cite{GB}, \cite{BM2}. In \cite{BM2}, it was shown that minimizing
sequences for the corresponding minimization problem develop a novel
type of so-called "near-boundary" vortices, which approach the
boundary and have finite GL energy in the limit of small $\ve$ (due
to the ghost vortices, see Appendix A). However, such minimizing
sequences do not converge to actual minimizers \cite{BGR}. These
studies lead to the natural question of whether there exist true
solutions of (\ref{equation}, \ref{boundary}) with near-boundary
vortices. Unlike the minimizing sequences such solutions may model
observable states of a physical system (e.g., persistent currents
with vortices and superfluids between rotating cylinders \cite{F}).
The following theorem, which is the main result of this work,
provides the answer to this question.

%This is
%done by developing appropriate variational technique based on the
%Direct method of Calculus of Variations.

\begin{thm}[Existence of solutions with
vortices of problem (\ref{equation}, \ref{boundary})]
\label{semistiffsolutions} For any integer $M>0$, there exist at
least $M$ distinct stable solutions of (\ref{equation},
\ref{boundary}) with (nearboundary) vortices when $\ve<\ve_1$
\text{(}$\ve_1=\ve_1(M)>0$\text{)}.
%They are classified by
%their degrees on $\partial\Omega$, $\partial\omega$ and
%an additional integer parameter related to the approximate bulk
%$degree (defined below).
The vortices of these solutions are at distance $o(\ve)$ from the
boundary and have bounded GL energy in the limit $\ve\to 0$. The
solutions are stable in the sense that they are (local) minimizers
of (\ref{P1}) in $\C{J}$.
\end{thm}

%The stable solutions of (\ref{equation}), (\ref{boundary}) in Theorem
%\ref{semistiffsolutions} are local
%minimizers of the energy functional (\ref{P1}) in $\C{J}$.
%The stable solutions of (\ref{equation}), (\ref{boundary}) we obtain in this work
%are local minimizers of the energy functional (\ref{P1}) in $\C{J}$.
To construct local minimizers of (\ref{P1}) in $\C{J}$, we represent
$\C{J}$ as the union of subsets
$\C{J}^{(d)}_{pq}$ (defined in (\ref{setJpqd}) below), $\C{J}=\cup_{p,q,d\in\D{Z}}\C{J}_{pq}^{(d)}$,
and study the existence of global minimizers in $\C{J}^{(d)}_{pq}$.
Furthermore, we show that each minimizer lies in $\C{J}^{(d)}_{pq}$ with its open
neighborhood.
%that these global minimizers in $\C{J}^{(d)}_{pq}$
%are interior elements
%of $\C{J}^{(d)}_{pq}$
Therefore, the minimizers in $\C{J}^{(d)}_{pq}$ are distinct local minimizers in $\C{J}$.

Thus, the construction of solutions of (\ref{equation},
\ref{boundary}) is based on the study of the following constrained
minimization problem:
\begin{equation}
\label{constrainedminimization}
m_\ve(p,q,d):=\inf\{E_\ve(u);\,u\in \C{J}_{pq}^{(d)}\},
\end{equation}
where
\begin{equation}
\label{setJpqd}
\C{J}_{pq}^{(d)}=\{u\in\C{J}_{pq}; {\rm abdeg}(u)\in [d-1/2,d+1/2]\},
\end{equation}
$p$,$q$ and $d$ are given integers, and ${\rm abdeg}(u)$ is the {\it approximate bulk degree},
introduced as follows. Consider the boundary value problem
\begin{equation}
\label{capacitysolution}
\begin{cases}
\Delta V=0 \quad \text{in}\ A\\
V=1\quad \text{on}\ \partial\Omega\\
V=0\quad \text{on}\ \partial\omega.
\end{cases}
\end{equation}
%By means of the solution $V$ of (\ref{capacitysolution}) we
Introduce
%Given a function $V\in C^1(\bar A)$
%such that $V=1$ on $\partial \Omega$ and $V=0$ on $\partial \omega$, introduce
${\rm abdeg}(\ \cdot\ ):\ H^1(A;\D{R}^2)\to\D{R}$ by the formula
\begin{equation}
{\rm abdeg}(u)=
\frac{1}{2\pi}\int_{A}
u\times(\partial_{x_1} V\, \partial_{x_2} u-
 \partial_{x_2} V\, \partial_{x_1} u)
\,{\rm d}x,
\label{adeg}
\end{equation}
where $V$ solves (\ref{capacitysolution}).
%In Sec. ? we show that ${\rm adeg}(u,A)$ is in some sense
In the particular case where $A$ is a circular annulus, $A_{R_1R_2}=\{x;R_1<|x|<R_2\}$,
${\rm abdeg}(u)$ is expressed by
\begin{equation}
{\rm abdeg}(u)=\frac{1}{\log({R_2}/{R_1})}\int_{R_1}^{R_2}
\Bigl(\frac{1}{2 \pi}\int_{|x|=\xi} u\times\frac{\partial
u}{\partial \tau}{\rm d}s \Bigr) \frac{{\rm d}\xi}{\xi}.
\label{1.9'}
\end{equation}

For $S^1$-valued maps,
${\rm abdeg}(u)$ becomes integer valued  and representation
(\ref{1.9'}) clarifies its interpretation
%of the ${\rm abdeg}(u)$
as an average value of the standard degree. The definition
(\ref{adeg}) is motivated by the following intuitive consideration:
represent the standard degree over the boundary $\partial \Omega$
via a "bulk" integral over the area of $A$ for $S^1$-valued maps and
notice that if $E_\ve(u)\leq \Lambda$ for some finite $\Lambda$ and
sufficiently small $\ve$, then $u$ is "almost" $S^1$-valued.

It was observed in \cite{A} that for $S^1$-valued maps in an annulus
$A$, one can define the topological degree ${\rm deg}(u,A)$ that
classifies maps $u\in H^1(A;S^1)$ according to their $1$-homotopy
type \cite{White} ($1$-homotopy type is completely determined by the
degree of the restriction to a nontrivial contour). This definition
was relaxed in \cite{A} for maps that are not necessarily
$S^1$-valued by considering $u/|u|$ in a subdomain $A_u\subset A$,
which is obtained by removing neighborhoods of the boundary
$\partial A$ and zeros (vortices) of $u$. Definition (\ref{1.9'})
does not require the removal of vortices from $A$, and ${\rm
abdeg}(u)$ is obtained by a simple formula (unlike ${\rm deg}(u,A)$
in \cite{A}, where the domain of integration depends on $u$). Note
that in general ${\rm abdeg}(u)$ is not an integer. The most
important fact for our consideration is that ${\rm abdeg}(u)$ is
continuous with respect to weak $H^1$-convergence, unlike the
standard degree in (\ref{def_degree}) (this issue for ${\rm
deg}(u,A)$ was not addressed in \cite{A}).

%a definition of a global degree
%${\rm deg}(u,A) \in\D{Z}$
%of a map on an annulus was introduced. While this
%definition has some common features with (\ref{1.9'}) there are
%essential differences which are discussed in Section
%\ref{section3}. In particular, ${\rm deg}(u,A)$ is defined by
%integration over a subdomain $A_u\subset A$ which depends on $u$
%in a quite complicated way.

%A similar formula is used in \cite{A} to define a generalized
%degree over an annulus. The relation the latter notion and the
%notion of the approximate bulk degree is discussed in Section
%\ref{section3}.

The minimization in problem (\ref{constrainedminimization})
is taken over $\C{J}_{pq}^{(d)}$, which is
not an open set, and therefore minimizers of
(\ref{constrainedminimization}) (if they exist) are not necessarily
local minimizers of (\ref{P1}) in $\C{J}$. Indeed, while
$\C{J}_{pq}$ is an open subset of $\C{J}$
(hereafter we assume the topology and convergence in
$\C{J}$ to be the strong $H^1$ unless otherwise is specified),
the constraint ${\rm abdeg}(u)\in[d-1/2,d+1/2]$ defines a closed set with respect
to both strong and weak $H^1$-convergences.
However,  if we further consider a subset of  maps
with bounded energy and choose $\ve$ small enough, then the constraint
${\rm abdeg}(u)\in[d-1/2,d+1/2]$ becomes open
%%when $\ve$ is small  and GL energy of
%this last possibility is excluded for
%low energy solutions of
%
thanks to the following proposition:
\begin{prop}
\label{openconstraint} Fix $\Lambda>0$. There exists
$\ve_0=\ve_0(\Lambda)>0$ such that if $0<\ve<\ve_0$,
then for any integer $d$ and any $u\in H^1(A;\D{R}^2)$ satisfying
$E_\ve(u)\leq\Lambda$ the closed constraint ${\rm abdeg}(u)\in[d-1/2,d+1/2]$ is equivalent
to an open one that is,
\begin{equation}
d-1/2\leq{\rm abdeg}(u)\leq d+1/2 \ \Longleftrightarrow\ d-1/2
<{\rm abdeg}(u)< d+1/2.
\end{equation}
\end{prop}

The following theorem is the main tool in proving existence of local minimizers.

\begin{thm}[Existence of minimizers of the constrained problem]
\label{mainth}
For any integers $p$, $q$ and $d>0$ ($d<0$)
with $d\geq\max\{p,q\}$ ($d\leq\min\{p,q\}$)
there exists $\ve_1=\ve_1(p,q,d)>0$
such that
the infimum in (\ref{constrainedminimization}) is
always attained,
%at a $u_\ve\in {\C J}_{pq}^{(d)}$
%$\Lambda>2(\pi d)^2/{\rm cap}(A)+\pi(|d-p|+|d-q|)$ and
when $\ve<\ve_1$. Moreover
\begin{equation}
m_\ve(p,q,d)\leq I_0(d,A)+\pi(|d-p|+|d-q|),
\label{add}
\end{equation}
where
\begin{equation}
I_0(d,A)=\min\left\{\frac{1}{2}\int_A|\nabla u|^2{\rm d}x,\ u\in H^1(A,S^1)\cap
\C{J}_{dd}\right\}.
\label{S^1minimization}
\end{equation}
The value $I_0(d,A)$ is expressed by $I_0(d,A)=2(\pi d)^2/{\rm cap}(A)$ via
the $H^1$-capacity ${\rm cap}(A)$ of the domain $A$.
\end{thm}

The key difficulty is to establish the attainability of the
infimum in (\ref{constrainedminimization}), which is highly
nontrivial since the degree on $\partial\Omega$ and $\partial\omega$ is not
preserved in the weak $H^1$-limit \cite{BM2,BGR}.
We show solutions of (\ref{equation}, \ref{boundary}), which are minimizers
of (\ref{constrainedminimization}) (local minimizers of (\ref{P1}) in $\C{J}$)
with $p\not=d$ and any $q$ (or $q\not=d$ and any $p$) must have
vortices.
For fixed $\ve$, these vortices are located at a positive distance
from $\partial A$ and approach $\partial A$ as $\ve\to 0$.

Without loss of generality, throughout
this work
we always assume that $d>0$ (otherwise one can reverse the orientation of $\D{R}^2$).

Theorem \ref{semistiffsolutions} follows from  Theorem \ref{mainth} and Proposition
\ref{openconstraint}. The asymptotic behavior of the local minimizers is established in

\begin{thm}[Asymptotic behavior of minimizers and their energies]
\label{asymptoticth} Assume that the integers  $p$, $q$ and $d$
satisfy the assumptions of Theorem \ref{mainth}. Then as $\ve\to 0$
minimizers of (\ref{constrainedminimization}) converge weakly in
$H^1(A)$, up to a subsequence, to a harmonic map $u$ which minimizes
(\ref{S^1minimization}). Additionally,
\begin{equation}
\label{asymptoticforenergy1}
E_\ve(u_\ve)=I_0(d,A)
%\frac{1}{2}\int_A|\nabla u|^2{\rm d}x
+\pi(|d-p|+|d-q|)+o(1),\quad \text{as}\ \ve\to 0\text{, and}
\end{equation}
\begin{equation}
\label{asymptoticforenergy2}
E_\ve(u_\ve)=\frac{1}{2}\int_A|\nabla u_\ve|^2{\rm d}x+o(1),
\quad \text{as}\ \ve\to 0.
\end{equation}
\end{thm}

In particular, it follows  from (\ref{asymptoticforenergy1},
\ref{asymptoticforenergy2}) that there is no strong convergence of
minimizers of (\ref{constrainedminimization}) in $H^1(A)$ as
$\ve\to
0$ unless $p=q=d$.

Next, we summarize the distinct features of the GL
boundary value problem with semi-stiff boundary conditions. The
first interesting feature is the existence of solutions with a new
type of vortices called {\it near-boundary vortices}. Unlike the
inner vortices, whose energy blows up at the rate of $|\log\ve|$,
the energy of near-boundary vortices is bounded as $\ve \to 0$ and they are
located at a distance $o(\ve)$ from the boundary.

Secondly, the semi-stiff boundary conditions result in a lack of
compactness. Namely, as of now, the only way to find nonconstant
minimizers is by searching for minimizers in subsets
$\C{J}_{pq}^{(d)}\subset \C{J}$. These subsets, however, are not
weakly $H^1$-closed and therefore a weak $H^1$-limit of a
minimizing sequence $(u^{(k)})\in \C{J}_{pq}^{(d)}$ may not lie in
$\C{J}_{pq}^{(d)}$, but rather in $\C{J}_{p^\prime q^\prime}^{(d)}$
with $p^\prime\not=p$ or (and) $q^\prime\not=p$. Theorem
\ref{mainth} shows that if $d>0$ and $d\geq\max\{p,q\}$, then (for
small $\ve$) any weak $H^1$-limit of a {\it minimizing sequence}
$(u^{(k)})\subset \C{J}_{pq}^{(d)}$ always belongs to
$\C{J}_{pq}^{(d)}$, despite the lack of weak $H^1$-closeness of
$\C{J}_{pq}^{(d)}$. In contrast, if $d\geq 0$ and  $d<\max\{p,q\}$
we have

\begin{Con}
Let $d\geq 0$, $d<\max\{p,q\}$  (or $d\leq 0$, $d>\min\{p,q\}$) and
let $u$ be a weak limit of a minimizing sequence for problem
(\ref{constrainedminimization}) (such a minimizing sequence exists
and bound (\ref{add}) holds for any integer $p$, $q$, $d$, see
Appendix B). Then $u\not\in\C{J}_{pq}^{(d)}$ when $\ve$ is
sufficiently small.
\end{Con}

In the simplest case, when $d=0$ and either $p=1$ and $q=0$ or $p=0$ and $q=1$, this conjecture is demonstrated
by an argument quite similar to the nonexistence proof in  \cite{BM1} for simply connected domains
(see Sec.\ref{section2} below). A more interesting example, which supports the  above conjecture, follows from  the previously studied (global) minimization problem
$\tilde m_\ve=\inf\{E_\ve(u),\,u\in\C{J}_{11}\}$. It was shown
in \cite{BM1,BGR}  that if ${\rm cap}(A)\geq\pi$ (subcritical/critical cases),
then $\tilde m_\ve$ is always attained, whereas if ${\rm cap}(A)<\pi$ (supercritical case),
then $\tilde m_\ve$ is never attained for small $\ve$. One can see elements of
minimizing sequences lie in $\C{J}^{(1)}_{11}$ in  subcritical/critical cases
and in $\C{J}^{(0)}_{11}$ in the supercritical case. Moreover, the nonexistence of
minimizers in problem (\ref{constrainedminimization}) for $d=0$, $p=q=1$ and small $\ve$
holds for any doubly connected domain (with any capacity). (For ${\rm cap}(A)<\pi$ the
proof is presented in \cite{BGR}, this proof can be easily generalized for
${\rm cap}(A)\geq\pi$.)

We conclude the introduction by outlining the scheme of the proof of
Theorem \ref{mainth}, which employs a comparison argument. Fix an
integer $d>0$. First, we establish the existence of minimizers in
problem (\ref{constrainedminimization}) for $p=q=d$ by using the
so-called Price Lemma \cite{BM2} (see Lemma \ref{pricelemma} below),
the uniform lower energy bound from Lemma \ref{roughlowerbound} and
the upper bound from Lemma \ref{ochenprostayalemma}, which is
obtained by considering $S^1$-valued testing maps. We show these
minimizers (which belong to $\C{J}_{dd}^{(d)}$)
%in $\C{J}_{dd}^{(d)}$
are vortexless. Next, we argue by induction on the parameter $\ae(p,q)=|d-p|+|d-q|$.
This parameter is naturally associated with the number of vortices--
for example, for the above minimizers in $\C{J}_{dd}^{(d)}$, we have $\ae(d,d)=0$.
Given an integer $K\geq 0$,
we assume the existence of minimizers in problem (\ref{constrainedminimization})
for $p$, $q$ such that, $\ae(p,q)\leq K$ and $p\leq d$, $q\leq d$
(the induction hypothesis) and prove the existence of minimizers
for $p$, $q$ such that $\ae(p,q)=K+1$ and $p\leq d$, $q\leq d$. The first step
in the induction procedure (when $K=0$) is shown in Section \ref{section5}.
The key technical point there is to construct a testing map $v\in \C{J}^{(d)}_{d(d-1)}$
such that
\begin{equation}
\label{short1}
E_\ve(v)<E_\ve(u_0)+\pi,
\end{equation}
where $u_0$ is a minimizer of (\ref{P1}) in $\C{J}^{(d)}_{dd}$.
This map $v$ is constructed by using the minimizer $u_0$ and
M$\rm\ddot{o}$bius conformal maps (Blashke factor \cite{BerVoss})
on the unit disk with a prescribed single zero near the boundary.
Then, given a minimizing sequence
$(u^{(k)})\subset\C{J}_{d(d-1)}^{(d)}$ of problem
(\ref{constrainedminimization}) for $p=d$, $q=d-1$, we have, by
(\ref{price2}) from Lemma [\ref{pricelemma}] and (\ref{short1}),
\begin{equation}
\label{short2}
E_\ve(u)+\pi(|d-{\rm deg}(u,\partial\omega)|+|d-1-{\rm deg}(u,\partial\Omega)|)\leq
\lim_{k\to\infty}E_\ve(u^{(k)})<E_\ve(u_0)+\pi,
\end{equation}
where u is a weak $H^1$-limit of  $(u^{(k)})$ (possibly a
subsequence). Then we estimate the left hand side  of (\ref{short2})
by the lower energy bound from Lemma \ref{roughlowerbound} and the
right hand side of (\ref{short2}) by the upper bound from Lemma
\ref{ochenprostayalemma}, this yields
\begin{equation}
\label{short3}
I_0(d,A)+\pi(2|d-{\rm deg}(u,\partial\omega)|+|d-1-{\rm deg}(u,\partial\Omega)|+|d-{\rm deg}(u,\partial\Omega)|)<
I_0(d,A)+\frac{3}{2}\pi.
\end{equation}
This implies that
${\rm deg}(u,\partial\omega)=d$, and either
${\rm deg}(u,\partial\Omega)=d-1$ or ${\rm deg}(u,\partial\Omega)=d$. In
view of (\ref{short2}), the only possible case is actually
${\rm deg}(u,\partial\Omega)=d-1$ since, otherwise, $u\in\C{J}_{dd}^{(d)}$
and therefore $E_\ve(u)\geq E_\ve(u_0)$ which contradicts
(\ref{short2}). Thus $u\in\C{J}_{d(d-1)}^{(d)}$ and is a minimizer in $\C{J}_{d(d-1)}^{(d)}$. The proof of existence
of minimizers for $p=d-1$, $q=d$ is quite similar.
So we have shown
that the existence of minimizers for $\ae(p,q)=0$ implies the existence of minimizers for $\ae(p,q)=1$.
In the
general case, when passing from $\ae(p,q)\leq K$ to
$\ae(p,q)\leq K+1$ in problem (\ref{constrainedminimization}), we
use the same idea but it is technically much more involved. It
requires the asymptotic analysis as $\ve\to 0$ of minimizers
$u_{pq}$ of (\ref{constrainedminimization}) with $\ae(p,q)=K$,
which is carried out in Section \ref{section6}. Based on the
result of this asymptotic analysis, we construct testing maps
$v\in\C{J}^{(d)}_{p^\prime q^\prime}$ ($p^\prime=p$,
$q^\prime=q-1$ or $p^\prime=p$, $q^\prime=q-1$), such that
$E_\ve(v)<E_\ve(u_{pq})+\pi$.

%in Section \ref{section2} of this work)
%The existence of minimizers is shown by a comparison argument. The main technical
%step is to show how to construct from a given minimizer $u$ in $\C{J}_{pq}^{(d)}$ a testing
%map $v$ in $\C{J}_{(p-1)q}^{(d)}$ (and in $\C{J}_{p(q-1)}^{(d)}$) such
%that $E_\ve(v)-E_\ve(u)<\pi$. Then, by using so-called Price Lemma \cite{BM2}
%(see Lemma \ref{pricelemma} in Section \ref{section2} of this work) and uniform
%lower bounds provided by Lemma \ref{roughlowerbound} (see Section \ref{section4})
%we get that any minimizing sequence in $\C{J}_{(p-1)q}^{(d)}$
%(or in $\C{J}_{p(q-1)}^{(d)}$) converges to a minimizer in this set.
%In Section \ref{section5} we show first existence of a
%(vortexless) minimizer in $\C{J}_{dd}^{(d)}$, then apply the above
%described procedure to show the existence of minimizers in
%$\C{J}^{(d)}_{(d-1)d}$ and in $\C{J}^{(d)}_{d(d-1)}$. These minimizers have
%exactly one vortex. In the general case of $p$, $q$ as in Theorem
%\ref{mainth} the proof of the existence of minimizers is carried out
%through Sections \ref{section6}, \ref{section7} by induction.

\section{Preliminaries}
\label{section2}

Throughout the paper we use the following notations.
\begin{itemize}
\item The vectors $a=(a_1,a_2)$ are identified with complex numbers $a=a_1+i a_2$.
\item $a\cdot b$ stands for the scalar product
$a\cdot b=a_1 b_1+a_2 b_2=\frac{1}{2}(a\bar b+\bar a b)$.
\item $a\times b$ stands for the vector product
$a\times b=a_1 b_2-a_2 b_1=\frac{i}{2}(a\bar b-\bar a b)$.
\item The orientation of simple (without self intersecting) curves
in $\D{R}^2$ (in particular $\partial \omega$ and $\partial\Omega$)
is assumed counterclockwise. If $\C{L}$ is such a curve, $\tau$ stands for the unit
tangent vector pointing in the sense of the above mentioned orientation on $\C{L}$,
$\nu$ is the unit normal vector such that $(\nu,\tau)$ is direct.
\item If $h$ is a scalar function, then $\nabla^\bot h =(-\partial_{x_2}h,\partial_{x_1}h)$.
\end{itemize}

\subsection{Properties of solutions from problem (\ref{equation}, \ref{boundary})}

As shown in \cite{BM1} by a bootstrap argument, any solution $u\in
H^1(A;\D{R}^2)$ of problem (\ref{equation}, \ref{boundary}) is
sufficiently regular (e.g., $u\in C^2(\overline{A})$ if
$A$ has a $C^2$ boundary). By the maximum principle we also have
\begin{lem} \label{maximummodulus}
The function $\rho(x)=|u(x)|$ satisfies $\rho\leq 1$ in $A$.
\end{lem}
Locally, away from its zeros, $u$ can be written as $u=\rho {\rm e}^{i\phi}$ with
real-valued phase $\phi$.
%This
%representation cannot, in general, be extended to the whole $A$.
We also will frequently make use of
the current potential $h$ related to the solution $u$
of (\ref{equation}, \ref{boundary}) by
\begin{equation}
\begin{cases}
\nabla^\bot h=(u\times \partial_{x_1} u,u\times \partial_{x_2} u)\quad \text{in}\ A\\
h=1\quad \text{on}\ \partial{\Omega}.
\end{cases}
\label{system_for_r}
\end{equation}
Unlike the phase $\phi$, the function $h$ is defined globally on $A$, and
\begin{equation}
\label{hphirelation}
\nabla h=-\rho^2\nabla^\bot \phi\quad \text{when}\ \rho>0.
\end{equation}
The existence of the unique solution of system (\ref{system_for_r}) and its elementary properties are
established in the following
\begin{lem} There exists the unique solution $h$ of the system (\ref{system_for_r}) and
$h={\rm Const}$ on $\partial \omega$, moreover
\begin{equation}
\Delta h=2\partial_{x_1} u\times\partial_{x_2} u\quad \text{in}
\ A,
\label{svoistvo_r1}
\end{equation}
\begin{equation}
\div(\frac{1}{\rho^2}\nabla h)=0 \quad \text{when}\ \rho>0.
\label{svoistvo_r2}
\end{equation}
\end{lem}
{\bf Proof.} The vector field
$F=(u\times \partial_{x_1} u,u\times \partial_{x_2} u)$ is divergence free.
Indeed, since $u$ is a smooth solution of  (\ref{equation}), we have
${\rm div}F=u\times \Delta u=0$
in $A$.
It follows that for any simply connected domain $W\subset A$ there is a unique
(up to an additive constant) function $\Phi$ solving
$\nabla^\bot \Phi=F$ in $W$ (this is well known Poincar\'e's lemma) .
Such a local solution $\Phi$ can be extended to a (possibly multi-valued) solution on $A$. Thanks to the
fact that $u$ satisfies (\ref{boundary}) we have $\frac{\partial \Phi}{\partial \tau}=-F\cdot\nu=0$ on
$\partial A$, i.e. $\Phi$ takes constant values on every connected component of
the boundary. This means $\Phi$ is actually a single valued function. Then,
$h(x)=\Phi(x)-\Phi(\partial\Omega)+1$ is the unique solution of (\ref{system_for_r}).
%$\Phi$
%satisfying $\Phi=1$ on  $\partial A\cap\partial W$. It follows that solutions $\Phi$ on different subdomains
%$W$ coincide with the same globally defined function.
%is just the restriction of
%$$Now, cover $A$ by two simply connected domains $W_k\subset A$ such that
%$\partial \Omega\cap\partial W_k\not=0$ and $\partial \Omega\cap\partial W_1\cap\partial %W_{k^\prime}\not=0$ if
%$W_k\cap W_{k^\prime}\not=\emptyset$ . Let $\Phi_k$ solve $\frac{\partial \Phi_k}{\partial x_1}=-F_2$,
%${\partial \Phi_k}{\partial x_2}=F_1$ in $W_k$ under the addition condition that
%$\frac{\partial \Phi_k}{\partial \tau}=0$ Thus there is a solution that can be extended to a (possibly %multivalued) function $r_\ve$ on $A$.
%On the other hand $\frac{\partial r_\ve}{\partial \tau}=u_\ve\times\frac{\partial u_\ve}{\partial \nu}=0$ on
%$\partial \Omega$, therefore (after adding a suitable constant) $r_\ve=1$ on $\partial\Omega$ and %$r_\ve={\rm Const}$
%on $\partial \omega$. This forces $r_\ve$ be a single-valued function.

The verification
of (\ref{svoistvo_r1}) is straightforward;
(\ref{svoistvo_r2}) follows directly from (\ref{hphirelation}).
\hfill$\square$
%write $u$ as $u=\rho {\rm e}^{i\phi}$ (locally when $\rho>0$)
%then
%(\ref{svoistvo_r2}) is nothing but
%$\frac{\partial^2 \phi}{\partial x_1\partial x_2}
%=\frac{\partial^2 \phi}{\partial x_2\partial x_1}$ that is obviously true.

In what follows, we also use the following result, which is valid
for any solution of the GL equation (\ref{equation}) (not
necessarily satisfying (\ref{boundary})).

\begin{lem}[\cite{M1}]
\label{Lem_Miron}
Let $u$ be a solutions of the GL equation (\ref{equation})
such that $|u|\leq 1$, $E_\ve(u)\leq \Lambda$,
where $\Lambda$ is independent of $\ve$. Then
\begin{equation}
\label{est_1_miron}
1-|u(x)|^2\leq \frac{\ve^2 C}{{\rm dist}^2(x,\partial A)}
\end{equation}
and
\begin{equation}
\label{est_2_miron}
|D^k u(x)|\leq \frac{C_k}{{\rm dist}^k(x,\partial A)}.
\end{equation}
where $C$, $C_k$ are independent of $\ve$.
\end{lem}

\subsection{Minimization among maps of $\C{J}$ with prescribed degrees}

\label{section22}

Any minimizer of (\ref{P1}) over the set $\C{J}_{pq}$
with prescribed integer degrees $p$ and $q$ is clearly a solution of
(\ref{equation}), (\ref{boundary}). However, the existence of minimizers is a nontrivial problem. In \cite{BV,GB,BM1,BM2,BGR}
the minimization problem for the Ginzburg-Landau functional (\ref{P1})
in $\C{J}_{11}$ was considered. In the case when $A$ is a circular annulus,
it was observed in \cite{BV} that minimizers, if they exist, brake the symmetry when
the ratio of the outer and inner radii of the
annulus exceeds certain threshold. By contrast, in the case when this ratio
is sufficiently close to $1$,
the existence of a unique minimizer and its symmetry is shown in \cite{GB}.
The techniques in both  \cite{BV} and \cite{GB} relied on the circular symmetry of
the domain. A more general approach based on the Price Lemma
was proposed in \cite{BM1},\cite{BM2}.

\begin{lem}
\label{pricelemma}
\cite{BM2}
Let $(u^{(k)}\in \C{J}_{pq})$ be a sequence that converges
to $u$ weakly in $H^1(A,\D{R}^2)$. Then
$$
\liminf_k\frac{1}{2}\int_A |\nabla u^{(k)}|^2 dx\geq
\frac{1}{2}\int_A |\nabla u|^2 dx+
\pi(|p-\deg (u,\partial \omega)|+|q-\deg (u,\partial \Omega)|),
$$
or, equivalently (by Sobolev embeddings),
\begin{equation}
\label{price2}
\liminf_k E_\ve(u^{(k)})\geq E_\ve(u)+ \pi(|p-\deg
(u,\partial \omega)|+|q-\deg (u,\partial \Omega)|).
\end{equation}
\end{lem}

This result is also of prime importance in this work. With the help of
Lemma \ref{pricelemma}, it was shown in \cite{BM2} that the infimum of
(\ref{P1}) in $\C{J}_{11}$ is always attained when ${\rm cap}(A)\geq\pi$.
It was also conjectured in \cite{BM2} that
when ${\rm cap}(A)<\pi$ and $\ve$ is
sufficiently small the weak limit of any minimizing sequence is not in
the class of admissible maps, i.e. the global minimizer does not
exist. In \cite{BGR} this nonexistence
conjecture was proved by a contradiction argument based on
explicit energy bounds.

While the existence/nonexistence of minimizers in $\C{J}_{pq}$ for $p=q=1$ is
nontrivial and the answer depends on ${\rm cap}(A)$ and also on $\ve$, the case
$p=0$, $q=1$ ($p=1$, $q=0$) is simple. Arguing as in \cite{BM1} we can
show that $\inf\{E_\ve(u); u\in\C{J}_{01}\}$ is never attained. Really,
we have
$$
\frac{1}{2}\int_{A}|\nabla u|^2{\rm d}x\geq
\left|\int_{A}\partial_{x_1} u\times\partial_{x_2}u
{\rm d}x\right|=\pi|{\rm deg}(u,\partial\Omega)-{\rm deg}(u,\partial\omega)|=\pi
$$
whenever $u\in\C{J}_{01}$. On the other hand, by constructing explicit sequence
in the spirit of \cite{BV} (see also \cite{BM1}), we have
$\inf\{E_\ve(u); u\in\C{J}_{01}\}=\pi$. Thus, if there exists a minimizer $u\in \C{J}_{01}$
then $u\in H^1(A;S^1)$ and $u$ solves the GL equation (\ref{equation}). This is impossible
unless $u$ is a constant map, and then $u\not\in \C{J}_{01}$.

\section{Properties of the approximate bulk degree}
%${\rm adeg}(u)$}
\label{section3}

The degree of restriction of maps from $H^1(A,S^1)$ to any smooth
closed curve, in particular $\partial\Omega$, is preserved by the weak
$H^1$-convergence. This result follows from \cite{White}, or can
be shown directly by using integration by parts as below in (\ref{degandadeg})
(note that, for any $S^1$-valued map $u$,
${\rm deg}(u,\partial\Omega)={\rm deg}(u,\partial\omega)={\rm deg}(u,\C{L})$, where $\C{L}$
is an arbitrary smooth simple curve in $A$ enclosing $\omega$).
Thus we have the decomposition
%$1$-homotopy type of smooth maps $u:\ A\to S^1$ is preserved by
%the weak $H^1$-convergence. The $1$-homotopy type of a smooth
%mapping from $A$ to $S^1$ is completely determined by the
%topological degree of its restriction to a smooth simple curve
%$\C{L}$ (with counterclockwise orientation) in $A$ enclosing
%$\omega$. Then, by the density of smooth maps  in $H^1(A,S^1)$
%(see, e.g., \cite{B}), we get the decomposition
\begin{equation}
\label{1-homotopy}
H^1(A,S^1)=\bigcup_{d\in \D{Z}}\{u\in H^1(A,S^1),\ {\rm deg}(u,\partial\Omega)=d\}
\end{equation}
by disjoint sets, each of them being closed in weak $H^1$-topology
of $H^1(A,S^1)$.
%Alternatively,  one can show directly that the function ${\rm deg}(\ \cdot\ ,\C{L})$
%is continuous with respect to the weak $H^1$-convergence on $H^1(A,S^1)$.
%(see,e.g. \cite{}, also ).

Fix $\Lambda>0$. In this section we consider maps $u \in H^1(A, \mathbb{R}^2)$ in the level
set $E^\Lambda_\ve=\{u;E_\ve(u)\leq \Lambda\}$. We show that the
approximate bulk degree ${\rm abdeg}(u)$ classifies maps $u\in
E^\Lambda_\ve$ similarly to the above classification
(\ref{1-homotopy}) of $S^1$-valued maps. The basic properties of
${\rm abdeg}(u)$ we demonstrate are
\begin{itemize}
\item[a)]${\rm abdeg}(u,A)=\deg(u,\partial \Omega)$ if $u\in H^1(A,S^1)$,
\item[b)]$|{\rm abdeg}(u)-{\rm abdeg}(v)|
\leq \frac{2}{\pi}\|V\|_{C^1(A)}\Lambda^{1/2}\|u-v\|_{L^2(A)}$
if $u,v\in E^\Lambda_\ve$.
\end{itemize}
The first property follows directly from the definition (\ref{adeg}) of ${\rm abdeg}(u)$. Really, integrating by parts in (\ref{adeg}),
we get
\begin{equation}
\label{degandadeg}
{\rm abdeg}(u)=\frac{1}{2\pi}\int_{\partial\Omega}u\times\frac{\partial u}{\partial \tau}
{\rm d}s-
\frac{1}{\pi}\int_{A}\partial_{x_1}u\times\partial_{x_2}u\,V {\rm d}x
=\deg(u,\partial \Omega),
%\in \D{Z}
\end{equation}
for any $u\in H^1(A,S^1)$ ($\partial_{x_1}u\times\partial_{x_2}u=0$ a.e. in $A$ since $|u|=1$ a.e.). The property b) of ${\rm abdeg}(u)$ is proved in the following
\begin{lem} For any $u,v\in H^1(A;\D{R}^2)$ we have
\label{cont}
$$
|{\rm abdeg}(u)-{\rm abdeg}(v)|
\leq \frac{1}{\pi}\|V\|_{C^1(A)}((E_\ve(u))^{1/2}+(E_\ve(v))^{1/2})\|u-v\|_{L^2(A)}.
$$
\end{lem}
{\bf Proof.}
% of Lemma \ref{cont}.}
We have, integrating by parts,
\begin{eqnarray*}
&2\pi({\rm abdeg}(u)-{\rm abdeg}(v))
&=\int_{A}
(u-v)\times
(\partial_{x_2} u \,\partial_{x_1} V  -\partial_{x_1} u \,\partial_{x_2} V
){\rm d}x
\\&&\ +\int_{A}v\times
(\partial_{x_2} (u-v)\,\partial_{x_1} V - \partial_{x_1} (u-v)\,\partial_{x_2} V
){\rm d}x
\\&&=\int_{A}(u-v)\times
(\partial_{x_2} u \,\partial_{x_1} V
-\partial_{x_1} u\, \partial_{x_2} V){\rm d}x
\\&&\ +\int_{A}
(u-v)\times
(\partial_{x_2} v \,\partial_{x_1} V
- \partial_{x_1} v\,\partial_{x_2} V){\rm d}x.
\end{eqnarray*}
Then the statement of the lemma follows by the Cauchy-Schwartz inequality.\hfill $\square$

The main consequence of properties a) and b) of the function ${\rm abdeg}(u)$ is

\noindent
\begin{prop}
\label{svoistvoc}
${\rm abdeg}(u)$ is close to integers uniformly in $u\in E^\Lambda_\ve$ when $\ve$ is
sufficiently small, i.e.
\begin{itemize}
\item[c)]
$\displaystyle
%\begin{equation*}
\sup_{u\in E_\ve^\Lambda}{\rm dist}({\rm abdeg}(u),\D{Z})\to 0$
as $\ve\to 0$.
\end{itemize}
\end{prop}
Before proving this fact note that Proposition \ref{svoistvoc}
immediately implies  Proposition \ref{openconstraint} stated in the Introduction.

\noindent
{\bf Proof of Proposition \ref{svoistvoc}.}
According to (\ref{degandadeg}) and Lemma \ref{cont} we have
\begin{equation}
\label{proofofpropertyc}
\sup_{u\in E_\ve^\Lambda}{\rm dist}({\rm abdeg}(u),\D{Z})\leq
\sup_{u\in E_\ve^\Lambda}\inf_{v\in E^\Lambda_0}
|\text{abdeg}(u)-\text{abdeg}(v)|\leq
\frac{2}{\pi}||V||_{C^1}\Lambda^{\frac{1}{2}}\delta_\ve,
\end{equation}
where $\delta_\ve$ is the following (nonsymmetric) distance
\begin{equation}
\delta_\ve:=\sup_{u\in E_\ve^\Lambda}{\rm dist}_{L^2(A)}(u,E^\Lambda_0)
\label{delta}
\end{equation}
between $E_\ve^\Lambda$ and
$$
E^\Lambda_0=
\left\{
u\in H^1(A,S^1);\ E_0(u)=\frac{1}{2}\int_A|\nabla u|^2{\rm d}x\leq\Lambda
\right\}.
$$
Show now that $\delta_\ve\to 0$ as $\ve\to 0$. In view of (\ref{proofofpropertyc})
this yields the desired result.
%implies property c) of ${\rm adeg}(u)$.
Assume by contradiction that $\delta_{\ve_k}\geq c>0$
for a sequence $\ve_k\to 0$. By virtue of the Sobolev embeddings
the supremum in (\ref{delta}) is attained on
$E_\ve^\Lambda$, i.e. $\delta_\ve={\rm
dist}_{L^2(A)}(u_\ve,E^\Lambda_0)$, where
$E_\ve(u_\ve)\leq\Lambda$. We can extract a subsequence of
$(u_{\ve_k})$, still denoted $(u_{\ve_k})$, that converges to a
map $u$ weakly in $H^1(A,\D{R}^2)$. Thanks to Sobolev embeddings
$u_{\ve_k}\to u$ strongly in $L^2(A)$ and we have $u\in
H^1(A,S^1)$, since $\int_{A}(|u_\ve|^2-1)^2dx\leq 4\Lambda\ve^2$ .
Besides $E_0(u)\leq \Lambda$, by the lower weak semicontinuity of the Dirichlet integral.
Thus $u\in E^\Lambda_0$ and
$\delta_{\ve_k}\leq\|u-u_{\ve_k}\|_{L^2(A)}\to 0$.\hfill $\square$

The following result illustrates the relation of ${\rm abdeg}(u)$,
which is not necessarily an integer with the standard notion of
degree on a curve. It provides a simple criterion whether the
constraint ${\rm abdeg}(u)\in[d-1/2,d+1/2]$ in
(\ref{constrainedminimization}) is satisfied in a particular case
when $u$ is a solution of equation (\ref{equation}).
% on $A$ the ${\rm abdeg}(u)\approx{\rm deg}(u,\C{L})$ on
%a contour in $\C{L}$

\begin{lem}
\label{sensadeg} Let $\C{L}=\{x\in A; V(x)=1/2\}$ denote the $1/2$
level set of $V$, where $V$ is the solution of
(\ref{capacitysolution}). ($\C{L}$ is a smooth curve enclosing $\omega$.)
%Choose the counterclockwise orientation on $\C{L}$.
Then if a solution $u$ of the GL equation (\ref{equation})
satisfies $|u|\leq 1$ in $A$ and $E_\ve(u)\leq\Lambda$, then {\rm
(i)} $|u|\geq 1/2$ on $\C{L}$ and {\rm (ii)} we have
\begin{equation}
{\rm abdeg}(u)\in[d-1/2,d+1/2]\ \Longleftrightarrow\ {\rm deg}(\frac{u}{|u|},\C{L})=d,
\end{equation}
when $\ve\leq\ve_1$, where $\ve_1=\ve_1(\Lambda)>0$ does not depend on $u$.
\end{lem}

{\bf Proof.} Consider the domain
\begin{equation}
\label{domainAdelta}
A_{\delta^\prime}=\{x\in A;\ \delta^\prime<V(x)<1-\delta^\prime\},
\end{equation}
 where
 $0<\delta^\prime<1/2$. It follows from Lemma \ref{Lem_Miron} that $u$ satisfies
 \begin{equation}
 \label{above1/2}
 |u|\geq 1/2\quad\text{on}\  A_{\delta^\prime},
 \end{equation}
when $\ve<\ve^\prime_1$
($\ve^\prime_1=\ve^\prime_1(\delta^\prime)>0$). This proves (i).
We can now write $u=\rho{\rm e}^{i\psi}$ ($\rho=|u|>1/2$) on
$A_{\delta^\prime}$. Find an extension of $\psi$ onto the whole
domain $A$. To this end we pass to a conformal image of $A$.
%$\C{O}$

It is well known (see,e.g. \cite{complanalysis}) that there is a
conformal mapping $\C{G}$
of $A$ onto the annulus $\C{O}$ with $R={\rm exp}(\pi / {\rm cap}(A))$ and $1/R$
as the outer and inner radii, correspondingly.
Moreover, $\C{G}$ is explicitly given by
$\C{G}={\rm exp}(\frac{2\pi}{{\rm cap}(A)}((V-1/2+i\Psi))$, where $\Psi$ is a
(multivalued) harmonic conjugate of $V$.  $\C{G}$ maps $A_{\delta^\prime}$ onto the
annulus $\C{G}(A_{\delta^\prime})\subset \C{O}$  whose outer and inner radii
are $R^\prime={\rm exp}(\frac{2\pi}{{\rm cap}(A)}(1/2-\delta^\prime))$
and $1/R^\prime$, correspondingly.

Now consider $\hat\psi(x)=\psi(G^{-1}(x))$ on $\C{G}(A_{\delta^\prime})$.
We can extend $\hat\psi$ to the whole $\C{O}$ by reflections
$\hat\psi(x):=\hat\psi(x(R^\prime)^2/|x|^2)$ when
$|x|\geq R^\prime$  and $\hat\psi(x):=\hat\psi(x/(R^\prime|x|)^2)$ when
$|x|\leq 1/R^\prime$, so that
\begin{equation}
\label{etensionproperties}
\int_{\C{O}}|\nabla\hat\psi|^2{\rm d}x\leq
\int_{\C{G}(A_{\delta^\prime})}|\nabla\hat\psi|^2{\rm d}x+
\int_{\C{O}\setminus\C{G}(A_{\delta^\prime})}|\nabla\hat\psi|^2{\rm d}x\leq 2\int_{\C{G}(A_{\delta^\prime})}|\nabla\hat\psi|^2{\rm d}x,
\end{equation}
when $0<\delta^\prime<1/4$. Then (\ref{etensionproperties}), the
conformal invariance of the Dirichlet integral and (\ref{above1/2}) imply
\begin{multline}
\label{bound007}
\int_{\C{O}}|\nabla\hat\psi|^2{\rm d}x\leq2\int_{\C{G}(A_{\delta^\prime})}|\nabla\hat\psi|^2{\rm d}x=
2\int_{A_{\delta^\prime}}|\nabla\psi|^2{\rm d}x\\
\leq 8
\int_{A_{\delta^\prime}}\rho^2|\nabla\psi|^2{\rm d}x\leq 16E_\ve(u)\leq 16\Lambda.
\end{multline}
The desired extension of $\psi$ onto $A$ is now given by  $\tilde\psi(x)=\hat\psi(\C{G}(x))$. Using again the
conformal invariance of the Dirichlet integral, we see by (\ref{bound007})
that
$E_\ve({\rm e}^{i\tilde\psi})=\frac{1}{2}\int_A|\nabla\tilde\psi|^2{\rm d}x
=\frac{1}{2}\int_{\C{O}}|\nabla\hat\psi|^2{\rm d}x
\leq 8\Lambda$. Besides, since
$\tilde\psi=\psi$ on $A_{\delta^\prime}$ and $\rho=|u|\leq 1$, we get
\begin{multline}
\label{normofdifference}
\|u-{\rm e}^{i\tilde\psi}\|_{L^2(A)}^2=
\int_{A\setminus A_{\delta^\prime}}|u-{\rm e}^{i\tilde\psi}|^2{\rm d}x+
\int_{A_{\delta^\prime}}(\rho-1)^2{\rm d}x\\
\leq 4|A\setminus A_{\delta^\prime}|+
\int_{A_{\delta^\prime}}(\rho^2-1)^2{\rm d}x\leq 4(|A\setminus A_{\delta^\prime}|+\Lambda\ve^2).
\end{multline}
Then, by choosing small $\delta^\prime$ and
$\ve_1=\ve_1(\delta^\prime)>0$(
$\ve_1(\delta^\prime)<\ve_1^\prime$), in view of Lemma \ref{cont},
bounds $E_\ve({\rm e}^{i\tilde\psi})\leq 8\Lambda$, $E_\ve(u)\leq
\Lambda$ and   (\ref{normofdifference}) we have $|{\rm
abdeg}(u)-{\rm abdeg}({\rm e}^{i\tilde\psi})|<1/2$, when
$\ve<\ve_1$. But ${\rm abdeg}({\rm e}^{i\tilde\psi})={\rm
deg}({\rm e}^{i\tilde\psi},\partial\Omega)= {\rm deg}({\rm
e}^{i\tilde\psi},\C{L})={\rm deg}(\frac{u}{|u|},\C{L})$, due to
(\ref{degandadeg}). Therefore if ${\rm
deg}(\frac{u}{|u|},\C{L})=d$ then $u\in E_\ve^{\Lambda,d}$, and
vice versa. Thus (ii) is proved. \hfill $\square$

\section{Minimization among $S^1$-valued maps. Some upper and lower
bounds for problem (\ref{constrainedminimization})}
\label{section4}

Consider the minimization problem
\begin{equation}
\label{S^1minimization'}
I_0(d,A^\prime):=\inf\{E_0(u);\,u\in H^1(A;S^1),\, {\rm deg}(u,\partial\omega^\prime)=
{\rm deg}(u,\partial\Omega^\prime)=d\},
\end{equation}
where $E_0(u)=\int_{A^\prime}|\nabla u|^2{\rm d}x$,
$A^\prime=\Omega^\prime\setminus\overline{\omega^\prime}$, and
$\omega^\prime$, $\Omega^\prime$ are smooth bounded simply connected domains in $\D{R}^2$,
such that $\overline{\omega^\prime}\subset \Omega^\prime$.
%The problem
%(\ref{H^1minimization'})
This problem is a particular case of the minimization problem considered in
\cite{BBH}(Chapter I).
%of the Dirichlet
%energy in the class of testing maps $u\in H^1(A,S^1)$
%having ${\rm deg}(u,\C{L})=d$ (cf. \ref{1-homotopy}). Since $u\in H^1(A,S^1)$ we have
%$$
%{\rm deg}(u,\C{L})={\rm deg}(u,\partial\omega)={\rm deg}(u,\partial\Omega)={\rm adeg}(u),
%$$
%i.e. instead of requiring ${\rm deg}(u,\C{L})=d$ we can equivalently redefine the class of
%testing maps as $\{u\in H^1(A,S^1),\ {\rm deg}(u,\partial\Omega)={\rm %deg}(u,\partial\omega)=d\}$.
%That is precisely in this setting the minimization problem is considered in \cite{BBH} %(Chapter I).
%The results we present below
%provide an (optimal) bound for (\ref{constrainedminimization}) in the case when
%$p=q=d$.
\begin{prop}\cite{BBH} There exists a unique (up to multiplication on constants with unit modulus)
solution $u$ of the minimization problem
(\ref{S^1minimization'}),
and $u$ is a regular harmonic map in $A^\prime$ (i.e.
$-\Delta u=u|\nabla u|^2$ in $A^\prime$)
satisfying $u\times\frac{\partial u}{\partial \nu}=0$ on $\partial A^\prime$.
\end{prop}
When $A^\prime=A$, then any minimizer $u$ of (\ref{S^1minimization'}) belongs to
$\C{J}_{dd}$.
By (\ref{degandadeg}) we also have ${\rm abdeg}(u)=d$. This yields the following
(optimal) bound for (\ref{constrainedminimization}), in the case when
$p=q=d$.
\begin{lem}
\label{ochenprostayalemma}
We have $m_\ve(d,d,d)\leq I_0(d,A)$ for
any $\ve>0$.
\end{lem}
It is shown in \cite{BBH} (Chapter I) that $I_0(d,A)$ can be expressed by
\begin{equation}
I_0(d,A)=\frac{1}{2}\int_A|\nabla h_0|^2{\rm d}x,
\end{equation}
where $h_0$ is the unique solution of the linear problem
\begin{equation}
\begin{cases}
\Delta h_0=0 \quad \text{in}\ A\\
h_0=1\quad \text{on}\ \partial \Omega,\quad h_0={\rm Const} \quad \text{on}\
\partial\omega\\
\displaystyle\int_{\partial \Omega} \frac{\partial h_0}{\partial
\nu}{\rm d}\sigma=2\pi d.
\end{cases}
\label{obschayaotsenka}
\end{equation}
$h_0$ and the solution $V$ of (\ref{capacitysolution}) are related via
$h_0=1+2\pi d(V-1)/{\rm cap}(A)$, where
${\rm cap}(A)$ stands for $H^1$-capacity of $A$ (see, e.g., \cite{M}). Thus
$I_0(d,A)=2(\pi d)^2/{\rm cap}(A)$, and this clearly holds for
any doubly connected domain
$A^\prime$ in place of $A$. Therefore we have
\begin{lem}
\label{capacitycont}
$I_0(d,A^\prime)$ depends continuously on ${\rm cap}(A^\prime)$.
\end{lem}
By using this simple result, we obtain the following lower bound for the
GL energy of solutions $u\in\C{J}$ of the equation
(\ref{equation}).
\begin{lem}
\label{roughlowerbound}
There is $\ve_2>0$ such that for any
solution $u\in \C{J}_{lm}$ of GL equation (\ref{equation})
satisfying $E_\ve(u)\leq\Lambda$, $abdeg(u)\in(d-1/2,d+1/2)$
we have
\begin{equation}
\label{roughboudfrombelow}
E_\ve(u)\geq I_0(d,A)-\frac{\pi}{2}+\pi(|d-l|+|d-m|),
\end{equation}
when $\ve<\ve_2$, and $\ve_2$ depends only on $\Lambda$.
% and $\hat\delta$.
\end{lem}

\noindent
{\bf Proof.} Due to the maximum principle $|u|\leq 1$ on $A$.
As in Lemma \ref{sensadeg} we consider the domain
$A_{\delta^\prime}$ that is defined by (\ref{domainAdelta}) and depends on a positive parameter $\delta^\prime<1/2$ to be chosen later.
Since $|u|\leq 1$ on $A$ we can apply Lemma \ref{Lem_Miron}  to get the bound
\begin{equation}
\label{bol'sheedinitsy} |u|\geq1-\ve\quad\text{in}\
A_{\delta^\prime},
\end{equation}
for $\ve<\ve^\prime_2$ ($\ve^\prime_2=\ve^\prime_2(\delta^\prime,\Lambda)>0$). Consider now the map
\begin{equation}
\label{newmap}
\tilde u= \frac{1}{1-\ve}\begin{cases}
u\quad\text{if}\  |u|<1-\ve\\
(1-\ve)\frac{u}{|u|}\ \text{otherwise}.
\end{cases}
\end{equation}
By (\ref{bol'sheedinitsy}) we have $|\tilde u|=1$ on $A_{\delta^\prime}$
and, according to Lemma \ref{sensadeg}, $\deg(\tilde u,\C{L})=d$ when $\ve<\min\{\ve_1,\ve_2^\prime\}$. Consequently,  the degree of $\tilde u$ on both connected
components of $\partial A_\delta^\prime$ equals $d$, so that $\frac{1}{2}\int_{A_{\delta^\prime}}|\nabla \tilde u|^2{\rm d}x\geq I_0(d,A_{\delta^\prime})$ (cf. (\ref{S^1minimization'})).
Therefore, by using the obvious pointwise inequality
$|\nabla \tilde u|^2\geq 2|\partial_{x_1}\tilde u\times\partial_{x_2}u|$ and the
integration by parts we get
\begin{multline}
\label{longbound}
\frac{1}{2}\int_{A}|\nabla \tilde u|^2{\rm d}x
\geq\frac{1}{2}\int_{A_{\delta^\prime}}|\nabla \tilde u|^2{\rm d}x
\\
+\sum_{k=1,2}\bigl|
\int_{A^{(k)}_{\delta^\prime}}
\partial_{x_1}\tilde u\times\partial_{x_2}\tilde u\,{\rm d}x\bigr|\geq I_0(d,A_{\delta^\prime})+\pi(|d-l|+|d-m|),
\end{multline}
where ${A^{(k)}_{\delta^\prime}}$ and $k=1,2$ are respectively the outer and
the inner connected components of $A\setminus A_{\delta^\prime}$. On the other
hand, it follows from (\ref{newmap}) that $|u|\leq|\tilde u|\leq1$.
Therefore,
%(\ref{newmap}) implies
\begin{equation}
\label{trickybound}
E_\ve(\tilde u)\leq\frac{1}{(1-\ve)^2}E_\ve(u).
\end{equation}
Bounds (\ref{longbound}) and (\ref{trickybound}) yield
(\ref{roughboudfrombelow}) when $\delta^\prime$ is such that
$I_0(d,A_{\delta^\prime})\geq I_0(d,A)-\pi/4$ (cf. Lemma
\ref{capacitycont}) and $\ve$ is sufficiently small.\hfill $\square$

%will be used in Sec. ? when deriving lower bounds for $E_\ve(u)$
%under the condition ${\rm adeg}(u)\in[d-1/2,d+1/2]$.

%{\bf Remark} It is clear that any minimizer $u$ of  (\ref{S^1minimization})
%belongs to ${\cal J}_{dd}^{(d)}(\ve,\Lambda)$, when $\Lambda>I_0(d,A)$.
%Let us also show that ${\cal J}_{dd}^{(d)}(\ve,\Lambda)\not=\emptyset$
%if $\Lambda=I_0(d,A))$.
%Indeed, let $u$ be a minimizer of  (\ref{S^1minimization}). Consider a map $v_\ve$ %minimizing
%$E_\ve(v)$ among $v\in H^1(A)$ with the Dirichlet boundary
%conditions $v=u$ on $\partial A$. According to \cite{BBH} and \cite{Struwe}
%$E_\ve(v_\ve)<I_0(d,A)$ for any $\ve>0$, and $v_\ve\to u$ strongly
%in $H^1(A)$ as $\ve\to 0$. It follows that
% $v_{\ve^\prime}\in{\cal J}_{dd}^{(d)}(\ve,\Lambda)$ when
% $\ve^\prime$ is sufficiently small.
%
%{\bf Remark ?.} In the case when  $A^\prime\subset A$ is a doubly connected
%subdomain of $A$, the same notations $I_0(d,A^\prime)$.

%\newpage

\section{From a vortexless minimizer to one with a single vortex}
\label{section5}
%{Vortexless and one vortex local minimizers}
%(global minimizers in ${\C J}_{dd}^{(d)}(\ve,\Lambda)$ and ${\C J}_{d(d-1)}^{(d)}(\ve,\Lambda)$)}

The main Theorem \ref{mainth} is proved by induction on
the ``number of vortices" in minimizers. More precisely,
given integer $d>0$, we show the existence of
minimizers of (\ref{constrainedminimization}) for $p=q=d$,
then pass to $p=d-1$, $q=d$ and $p=d$, $q=d-1$, e.t.c.
%
%More precisely, as
%we already shown in Section ? minimizing sequences of $E_\ve(u)$ in
%${\C J}_{pq}^{(d)}(\ve,\Lambda)$ may weakly converge to a
%map $u\in{\C J}_{p^\prime q^\prime}^{(d^\prime)}(\ve,\Lambda)$
%and
%
%we consider minimization  problem for $E_\ve(u)$ in ${\C J}_{pq}^{(d)}(\ve,\Lambda)$
%, the energy level $\Lambda$
%starting by proving existence i
The key point of the proof is the induction step, when the degree
changes by one on $\partial\omega$ or $\partial\Omega$. This change
results in the rise of an additional vortex in a minimizer. For
arbitrary $p$ and $q$, satisfying the conditions of Theorem
\ref{mainth}, this step is quite technical. This is why we consider
here a particular case of transition from $p=q=d$ (no vortices) to
$p=d$, $q=d-1$ (one vortex). The transition from $p=q=d$ to $p=d-1$
and $q=d$ is quite similar.
%We show that for sufficiently small $\ve$
%By using Price Lemma and continuity of ${\rm adeg}(u)$ with respect to
%weak $H^1$-convergence we establish
We first establish

%This number is discribed by the parameter
%$\ae(p,q)=|d-p|+|d-q|$, where $d>0$ is fixed, and $p\leq d$, $q\leq d$. The proof can be %divided into three parts. Part 1 (the induction base).
%establishing existence of minimizers of (\ref{constrainedminimization})
%for $p=q=d$. Part 2, induction step from
%$\ae(p,q)=k$ to $\ae(p,q)=k+1$. For general $k$ the proof is quite
%technical that is why we
%first present it the simplest case $k=0$, that is when
%passing from $p=q=d$ to $p=d-1$,$q=d$ and $p=d$, $q=d-1$.
%The key point of the proof is the induction step
%when degree change by one on $\partial\omega$ or $\partial\Omega$.
%We show that this transition results in a rise of an additional vortex in
%a minimizer.

%Part 2 (the first induction step) passing form
%$p=q=d$, then pass to $p=d-1$,$q=d$ and $p=d$, $q=d-1$, i.e.
%from $\ae(p,q)=0$ to $\ae(p,q)=1$. Part 3, induction step from
%$\ae(p,q)=k$ to $\ae(p,q)=k+1$.

\begin{lem}
\label{firststep} Given an integer $d>0$.
For sufficiently small $\ve$, $\ve\leq\ve_3$ with $\ve_3>0$,
the infimum $m_\ve(d,d,d)$ in (\ref{constrainedminimization}) is always attained,
and $m_\ve(d,d,d)\leq I_0(d,A)$.
%\frac{2(\pi d)^2}{{\rm cap}(A)}$.
%Then
%${\C J}_{dd}^{(d)}(\ve,\Lambda)\not=\emptyset$ and for sufficiently small
%$\ve$, $\ve\leq\ve_?$ with $\ve_?>0$,
% ${\rm inf}\{E_\ve(u);\, u\in{\C J}_{dd}^{(d)}(\ve,\Lambda)\}$ is always attained.
\end{lem}
\noindent
{\bf Proof.} Let $u$ be a weak $H^1$-limit of a minimizing sequence $(u^{(k)})$.
Since any minimizer $v$ of problem (\ref{S^1minimization}) is an
admissible testing map for problem (\ref{constrainedminimization}),
such a minimizing sequence exists and by using Price Lemma we obtain
%Assume that $\Lambda>\frac{2(\pi d)^2}{{\rm cap}(A)}$ then
%${\cal J}^{(d)}_{dd}(\ve,\Lambda)\not=\emptyset$, due to
%the above assumption
%on $\Lambda$ and
%results of Sec. 3.
%Let $u$ be a weak $H^1$-limit of a minimizing sequence $(u^{(k)})$.
%Then by the Price Lemma
%and Sobolev embedding theorems
%lower semicontinuity of $E_\ve(u)$
%we have
\begin{equation}
\label{contrbound}
E_\ve(u)+
%\pi(|{\rm deg}(u_\ve,\partial\Omega)-d|+|{\rm deg}(u_\ve,\partial\omega)-d|)
%\pi\ae({\rm deg}(u_\ve,\partial\Omega),{\rm deg}(u_\ve,\partial\omega))
\pi(|l-d| +|m-d|)
\leq\liminf_{k\to\infty}E_\ve(u^{(k)})\leq\frac{1}{2}\int_{A}|\nabla v|^2dx= I_0(d,A),
\end{equation}
where $l={\rm deg}(u,\partial\omega))$, $m={\rm deg}(u,\partial\Omega)$.
%(since $v\in H^1(A,S^1)$ and therefore
%$d={\rm deg}(v,\partial\Omega)=
%{\rm deg}(v,\partial\omega)={\rm adeg}(v,A)$, we have
%$v\in{\cal J}^{(d)}_{dd}(\ve,\Lambda)$).
Due to Proposition \ref{openconstraint} we have ${\rm abdeg}(u)\in (d-1/2,d+1/2)$
when $\ve<\ve_0$, therefore the first variation of
(\ref{P1}) at $u$ vanishes, i.e. $u$ is a solution of equation (\ref{equation}).
Indeed, thanks to Lemma  \ref{cont}, we have for any $w\in H^1_0(A;\D{R}^2)$ with sufficiently small $H^1$-norm, $u^{(k)}+w$ is an admissible testing map when $k$ is large, hence
$E_\ve(u+w)-E_\ve(u)=\lim_{k\to\infty}(E_\ve(u^{(k)}+w)-E_\ve(u^{(k)})\geq 0$
(where $\lim_{k\to\infty}$ denotes any partial limit), and we are done. Now, since $u$ is a solution of  (\ref{equation}), we can apply Lemma \ref{roughlowerbound}. That is we substitute
(\ref{roughboudfrombelow}) in (\ref{contrbound}), this leads to
\begin{equation}
|d-l| +|d-m|\leq \frac{1}{4},
\end{equation}
when $\ve<\min\{\ve_0,\ve_2\}$, i.e. $l=m=d$ (since $l$, $m$ and $d$ are integers).
Thus the infimum in (\ref{constrainedminimization}) for $p=q=d$ is always attained
when $\ve$ is sufficiently small. Lemma is proved.\hfill $\square$

Next, we perform the transition from the minimization
problem (\ref{constrainedminimization}) for $p=q=d$ to that for $p=d$, $q=d-1$
and show that $m_\ve(d,d-1,d)$ is always
attained when $\ve$ is sufficiently small. This is done by a comparison
argument of $m_\ve(d,d-1,d)$ with the
energy $E_\ve(u)$ of a minimizer $u$ of
(\ref{constrainedminimization}) for $p=q=d$.
We first describe the properties of such a minimizer.

%The proof of this Lemma is postponed till the end of this Section.

%Now we focus on the rise of vortex issue.

In Section 5, it is shown that for small $\ve$ any minimizer
$u$ of (\ref{constrainedminimization}) for $p=q=d$ is
vortexless (see Remark \ref{remvortexless}), i.e.
$u=\rho{\rm e}^{i\psi}$ with smooth $\rho>0$ and
$\psi:A\to \D{R}\setminus 2\pi d\D{Z}$ (torus).
It follows that we can write
$u=\rho{\rm e}^{id\theta}$,
the ${\rm e}^{i\theta}$ and $\nabla \theta$ being
smooth maps globally on $A$.
Then the boundary value problem
(\ref{equation})-(\ref{boundary}) rewritten in
terms of $\rho$ and $\theta$ is
%by the regularity results \cite{BerlMironPreprint}
\begin{equation}
\label{systempolarphase}
\begin{cases}
\displaystyle{\rm div}(\rho^2\nabla \theta) =0\quad \text{in}\ A\\
\displaystyle\frac{\partial\theta}{\partial\nu}=0\quad \text{on}\ \partial A,
\end{cases}
\end{equation}
\begin{equation}
\begin{cases}
\label{systempolarmodulus}
\displaystyle
-\Delta\rho+d^2|\nabla \theta|^2\rho+\frac{1}{\ve^2}\rho(\rho^2-1)=0
\quad \text{in}\ A\\
\displaystyle
\rho=1\quad \text{on}\ \partial A.
\end{cases}
\end{equation}
%As we will see minimizers in ${\C J}_{d(d-1)}^{(d)}(\ve,\Lambda)$ (if exist)
%have a vortex, unlike $u$ . It is located near boundary by contrast to ...
We also have, in view of (\ref{hphirelation}),
%\begin{eqnarray}
%\nabla h=-d\rho^2{\nabla}^\perp\theta &&
%\text{in}\ A \label{system_r1}\\
%h=1&& \text{on}\ \partial\Omega \label{system_r2},
%\end{eqnarray}
\begin{equation}
\nabla h=-d\rho^2{\nabla}^\perp\theta
\quad \text{in}\ A,
\label{system_r1}
%h=1&& \text{on}\ \partial\Omega \label{system_r2},
\end{equation}
where $h$ is the solution of (\ref{system_for_r}). It follows that $(h,\theta)$
defines orthogonal local coordinates in a neighborhood of $\partial\Omega$, thus straightening
out the boundary.
Really, it is straightforward
to verify that
$$
1-h(\partial\omega)=\frac{1}{{\rm cap}(A)}\int_A\nabla h\cdot \nabla V {\rm d}x=
2\pi{\rm abdeg}(u)/{\rm cap}(A),
$$
while ${\rm abdeg(u)}\geq(d-1/2)>0$. Then by applying the maximum
principle to (\ref{svoistvo_r2}) we get,
$1>h(x)>h(\partial\omega)$ in $A$. This, in turn, implies, by Hopf's
boundary lemma, that  $\frac{\partial h}{\partial\nu}>0$ on
$\partial\Omega$; i.e. the map
$(h,\theta):A\to\D{R}\times\D{R}\setminus2\pi\D{Z}$ can be
extended to a $C^1$-diffeomorphism  of a one sided neighborhood of
$\partial\Omega$  onto its image. Thus, there are some $\delta>0$
and a domain $G_\delta\subset A$ such that
$$
x\in G_\delta \to(h,\theta)\in \Pi_\delta=(1-\delta,1)\times\D{R}\setminus2\pi\D{Z}
$$
is a one-to-one correspondence, which extends to a $C^1$-diffeomorphsm of
$\overline{G}_\delta$ onto $[1-\delta,1]\times\D{R}\setminus2\pi\D{Z}$.

The following Proposition
is crucial for existence of minimizers of (\ref{constrainedminimization})
for $p=d$, $q=d-1$.
In particular, combined with Lemma \ref{ochenprostayalemma}
it provides an independent of $\ve$ bound for $m_\ve(d,d-1,d)$.

%We are going to show next that  ${\rm Inf}\{E_\ve(u);\, u\in{\C J}_{d(d
%-1)}^{(d)}(\ve,\Lambda)\}$
%is also attained, when $\Lambda>\frac{2(\pi d)^2}{{\rm cap}(A)}+\pi$. The following %Claim
%plays the crucial role in establishing existence of minimizers in
%${\C J}_{d(d-1)}^{(d)}(\ve,\Lambda)$ (in particular, it guaranties that this set is not %empty).

\noindent
\begin{prop}
\label{testingmap} Let $u=\rho {\rm e}^{id\theta}$, $\rho>0$, be a
minimizer of (\ref{constrainedminimization}) for $p=q=d$. Assume
that $\ve$ is so small that Proposition \ref{openconstraint}
holds with $\Lambda=I_0(d,A)$.
%and assume that $\rho_\ve:=|u_\ve|>0$ in $A$.
Then there is a testing map
$v\in{\C J}_{d(d-1)}$ such that ${\rm abdeg}(u)\in(d-1/2,d+1/2)$ and
\begin{equation}
\label{samayaglavnayaotsenka1}
E_\ve(v)-E_\ve(u)<\pi.
\end{equation}
\end{prop}

In Section \ref{section7}, the generalized version of Proposition \ref{testingmap} is used to show
existence of minimizers with several vortices.

\noindent
{\bf Proof of Proposition \ref{testingmap}.}
%For simplicity we suppress the subscript $\ve$.
We seek the testing map
$v$ in the form
\begin{equation}
\label{factorrepresent}
v=\rho w_t,
\end{equation}
with an unknown for the moment $w_t$. The following Lemma allows us to
compare the energy $E_\ve(u)$ of $u$ with that of $v$.

\begin{lem}
\label{comparisonlemma}
If $w\in H^1(G^\prime,\D{R}^2)$, $G^\prime\subset A$, is such that $|w|=1$ on $G^\prime$,
then
$$
\int_{G^\prime}(|\nabla(\rho w)|^2
+\frac{1}{2\ve^2}(|\rho w|^2-1)^2)dx=
\int_{G^\prime}(|\nabla u|^2
+\frac{1}{2\ve^2}(|u|^2-1)^2)dx+2L_\ve^{(d)}(w,G^\prime),
$$
where
\begin{equation}
\label{Lepsilon}
L_\ve^{(d)}(w,G^\prime)=\frac{1}{2}\int_{G^\prime}\rho^2|\nabla w|^2dx-
\frac{d^2}{2}\int_{G^\prime}
|\nabla \theta|^2 \rho^2 |w|^2dx
+\frac{1}{4\ve^2}\int_{G^\prime} \rho^4(|w|^2-1)^2dx
\end{equation}
\end{lem}

\noindent
This result is a variant of the factorization argument due to
\cite{CM}, its proof is presented in
the end of this section.

Note that if $G^\prime=G_\delta$ we can rewrite the functional
(\ref{Lepsilon})
by using local coordinates $(h,\theta)$ as (cf. (\ref{system_r1}))
\begin{multline}
L_\ve^{(d)}(w,G_\delta)=\frac{d}{2}\int_{\Pi_\delta}
|\partial_h w|^2\rho^2 {\rm d}h{\rm d}\theta\\
+\frac{1}{2d}\int_{\Pi_\delta}
(|\partial_\theta w|^2-d^2|w|^2){\rm d}h{\rm d}\theta
+\frac{1}{4\ve^2}\int_{\Pi_\delta} \rho^2(|w|^2-1)^2\frac{{\rm d}h{\rm d}\theta}
{d|\nabla\theta|^2}.
\end{multline}
Instead of dealing with  $L_\ve^{(d)}(w,G_\delta)$ we will
make use of the simplified functional
with a quadratic penalty term,
\begin{equation}
\label{Mlambda}
M_\lambda(w)=\frac{1}{2d}\int_{\Pi_\delta}
(d^2|\partial_h w|^2+|\partial_\theta w|^2){\rm d}h{\rm d}\theta
+\frac{1}{2d}\int_{\Pi_\delta}(\lambda|w-
{\rm e}^{i\theta}|^2-d^2|w|^2){\rm d}h{\rm d}\theta.
\end{equation}
This last functional admits the separation of variables.

Now consider the map $w_t$ that is given by $w_t={\rm e}^{id\theta}$ in
$A\setminus G_\delta$,
and continued to $G_\delta$ as a minimizer of the functional $M_\lambda(w)$,
where  $\lambda\geq 2d^2$, with the following prescribed boundary data:
\begin{equation}
\label{boundary11}
w_t={\rm e}^{id\theta} \C{F}_t({\rm e}^{i\theta})\quad\text{on}\ \partial\Omega,
\end{equation}
\begin{equation}
\label{boundary21}
w_t={\rm e}^{id\theta}\quad \text{on}\ \partial G_\delta\setminus \partial\Omega,
\end{equation}
where $\C{F}_t(z):=\C{C}_t(\bar z)$ (bar stands for the complex conjugate),
$\C{C}_t(z)=\frac{ z-(1-t)}{ z(1-t)-1}$
%the restriction to the boundary of
is the classical M$\rm\ddot{o}$bius conformal map from the unit
disk onto itself, $t<1$ is a positive parameter. Both parameters
$\lambda$ and $t$ will be determined later. Since ${\rm
deg}(\C{F}_t,S^1)=-1$ and ${\rm deg}({\rm
e}^{i\theta},\partial\Omega)=1$, the standard properties of the
topological degree implies that if $v$ is as in
(\ref{factorrepresent}), then
\begin{equation}
\label{nuzhnyedegree}
{\rm deg}(v,\partial\Omega)=d-1,\  {\rm deg}(v,\partial\omega)=d.
\end{equation}
The map $w_t$ is well defined now, because the functional $M_\lambda(w)$
with Dirichlet condition on the boundary has a unique minimizer for $\lambda\geq 2d^2$.
Moreover $|w_t|\leq 2$, since for if not then by
taking $\tilde w_t=\frac {w_t}{|w_t|}\min\{|w_t|,2\}$ in place of $w_t$ the first term in (\ref{Mlambda}) does not increase while the second one decreases, i.e.
$M_\lambda(\tilde w_t)<M_\lambda(w_t)$, it is a contradiction.

Note that under the following choice of $\lambda$,
$$
\lambda:=
\max\left\{\frac{9}{2\ve^2\inf_{G_\delta}{|\nabla \theta|^2}},
\,2d^2\right\}
$$
($|\nabla \theta|>0$ on the closure of $G_\delta$) we have
\begin{multline*}
\rho^2(|w_t|^2-1)^2\leq(|w_t|-1)^2(|w_t|+1)^2\leq|w_t-{\rm e}^{id\theta}|^2(|w_t|+1)^2\\
\leq
9|w_t-{\rm e}^{id\theta}|^2\leq 2\ve^2\lambda|\nabla\theta|^2|w_t-{\rm e}^{id\theta}|^2
\quad\text{in}\ G_\delta,
\end{multline*}
thanks to the bounds $|w_t|\leq 2$ and $\rho\leq 1$. It follows
that $L_\ve^{(d)}(w_t,G_\delta)\leq M_\lambda(w_t)$.  We get now, by
virtue of Lemma \ref{comparisonlemma}, that $v=\rho w_t$
%as in (\ref{factorrepresent})
satisfies
\begin{equation}
\label{bestotsenka}
E_\ve(v)\leq E_\ve(u)+M_\lambda(w_t).
\end{equation}

We can obtain a representation for
$M_\lambda(w_t)$ in separated variables. Namely, expanding $z^d\C{F}_t(z)$ on
$S^1$ as
$$
z^d\C{F}_t(z)=(1-t)z^d+
t(t-2)\sum_{k=0}^\infty (1-t)^kz^{d-k-1},
$$
we have
\begin{equation}
\label{ryaddlywt}
w_t=(1-tf_{-1}(h)){\rm e}^{id\theta} +
t(t-2)\sum^\infty_{k=0} (1-t)^k f_{k}(h){\rm e}^{-i(k-d+1)\theta},
\end{equation}
where $f_{k}(h)$ satisfy, according to (\ref{boundary11}, \ref{boundary21}),
\begin{equation}
\label{newboundcond}
f_{k}(1-\delta)=0,\quad f_k(1)=1.
\end{equation}
Substitute (\ref{ryaddlywt}) into (\ref{Mlambda}) to obtain
\begin{equation}
\label{mmm}
M_\lambda(w_t)=\frac{t^2\pi}{d} \Phi_{-1}(f_{-1})+
\frac{t^2\pi}{d}\sum_{k=0}^\infty (t-2)^2(1-t)^{2k}\Phi_k(f_k),
\end{equation}
where
\begin{equation}
\label{odnomernyi1}
\Phi_k(f_k)=\int_{1-\delta}^{1}
\bigl(d^2|f_{k}^{\prime}(h)|^2+
((k-d+1)^2+\lambda-d^2)|f_{k}(h)|^2\bigr) {\rm d}h.
\end{equation}
Minimizing (\ref{odnomernyi1}) under the conditions (\ref{newboundcond})
we get
\begin{equation}
\label{deffk}
f_k(h)=\frac{{\rm e}^{k_{+}(h-1)}}{1-{\rm e}^{(k_{-}-k_{+})\delta}}
+\frac{{\rm e}^{k_{-}(h-1)}}{1-{\rm e}^{(k_{+}-k_{-})\delta}},
\end{equation}
where $\displaystyle k_{\pm}=\pm\frac{1}{d}\sqrt{(k-d+1)^2+\lambda-d^2}$.
Therefore, we have
\begin{equation}
\label{asimptotika}
\Phi_k(f_k)= d(k-d+1)(1+\frac{\lambda-d^2}{2k^2}+O(1/k^3)),
\ \text{as}\ k\to\infty.
\end{equation}
Finally, using (\ref{asimptotika}) in (\ref{mmm}), we obtain
\begin{multline}
\label{finalotsenka11}
M_\lambda(w_t)\leq \pi((1-t)^2-1)^2\sum_{k=0}^\infty k(1-t)^{2k}
+2\pi t^2(\lambda-d^2)\sum_{k=1}^\infty \frac{(1-t)^{2k}}{k}+Ct^2\\
=\pi(1-2t-2t^2(\lambda-d^2)\log(1-(1-t)^2))+(C+\pi)t^2.
\end{multline}
Observe that the right hand side of (\ref{finalotsenka11}) is
strictly less than $\pi$ when $t>0$ is chosen sufficiently small.
By (\ref{bestotsenka}), for such $t$ the map $v=\rho w_t$ satisfies
(\ref{samayaglavnayaotsenka1}).

It remains only to show that ${\rm abdeg}(v)\in (d-1/2,d+1/2)$. To
this end note that by (\ref{ryaddlywt}) and  (\ref{deffk})
$w_t\to{\rm e}^{id\theta}$ pointwise in $G_\delta$ as $t\to 0$.
Therefore $\rho w_t\to \rho{\rm e}^{id\theta}(=u)$ weakly in
$H^1(A)$, so that ${\rm abdeg}(\rho w_t)\to{\rm abdeg}(u)$. On the
other hand, by Proposition \ref{openconstraint}, we know
$d-1/2< {\rm abdeg}(u)<d+1/2$. Thus, after possibly passing to a
smaller $t$, $v=\rho w_t$ satisfies the required property.
\hfill $\square$

Now we have that under the conditions of Proposition \ref{testingmap},
there exists a minimizing sequence $(u^{(k)})$ of admissible
testing maps in problem (\ref{constrainedminimization}) for $p=d$,
$q=d-1$ such that $\lim_{k\to\infty}
E_\ve(u^{(k)})<m_\ve(d,d,d)+\pi$ and $u^{(k)}$ weakly
$H^1$-converge to a map $u\in\C{J}$. Moreover, any minimizing
sequence has a subsequence with the same properties. Show that any
weak limit $u$ is also an admissible map. Let $l={\rm
deg}(u,\partial\omega)$, $m={\rm deg}(u,\partial\Omega)$. By
virtue of Lemma \ref{roughlowerbound} (in the same way us in Lemma
\ref{firststep} one shows that $u$ satisfies (\ref{equation})) and
Lemma \ref{pricelemma} we have
\begin{multline}
\label{posledsravnen}
I_0(d,A)-\frac{\pi}{2}+\pi(2|d-l|+|d-1-m|+|d-m|)\\
\leq E_\ve(u)
+\pi(|d-l|+|d-1-m|)<m_\ve(d,d,d)+\pi,
\end{multline}
since ${\rm abdeg}(u)=\lim_{k\to\infty} {\rm abdeg}(u^{(k)})\in
[d-1/2,d+1/2]$. Due to Lemma \ref{firststep} $m_\ve(d,d,d)\leq
I_0(d,A)$ so that (\ref{posledsravnen}) implies that $l=d$ and
either $m=d-1$ or $m=d$. In the last case, $u$ becomes an admissible map in problem
(\ref{constrainedminimization}) for $p=q=d$ and therefore
$E_\ve(u)\geq m_\ve(d,d,d)$, which contradicts the last
inequality in (\ref{posledsravnen}). Thus, $u\in\C{J}_{d(d-1)}$,
${\rm abdeg}(u)\in [d-1/2,d+1/2]$, i.e. $u$ is in the set of
admissible testing maps of problem (\ref{constrainedminimization})
for $p=d$, $q=d-1$.

\smallskip
\noindent
{\bf Proof of Lemma \ref{comparisonlemma}.} We have, by using (\ref{systempolarmodulus}),
\begin{multline*}
\int_{G_\delta}|\nabla(\rho w)|^2{\rm d}x= \int_{G_\delta}
(\rho^2|\nabla w|^2+\nabla\rho\cdot\nabla(\rho(|w|^2-1))+|\nabla\rho|^2){\rm d}x\\
=\int_{G_\delta}(\rho^2|\nabla w|^2+d^2\rho^2|\nabla
\theta|^2+\frac{1}{\ve^2}\rho^2(\rho^2-1)){\rm d}x\\
-\int_{G_\delta}(d^2\rho^2|\nabla \theta|^2
|w|^2+\frac{1}{\ve^2}\rho^2(\rho^2-1)|w|^2-|\nabla\rho|^2){\rm d}x.
\end{multline*}
Then simple algebraic manipulations give the required
result.\hfill $\square$

\section{Asymptotic behavior of local minimizers}

\label{section6}
In the previous section, we established the existence of minimizers
of (\ref{P1}) in $\C{J}^{(d)}_{dd}$
%problem (\ref{constrainedminimization}) for $p=q=d$
and demonstrated the first induction step of the proof of Theorem
\ref{mainth} that consists in transition from $p=q=d$ to $p=d$,
$q=d-1$ in (\ref{constrainedminimization}). (In fact, modulo the
assumption that  minimizers in $\C{J}^{(d)}_{dd}$ are vortexless,
we actually proved the existence of minimizers in $\C{J}^{(d)}_{(d-1)d}$ and
$\C{J}^{(d)}_{d(d-1)}$.) In order to show the induction step for
any integer $p\leq d$ and $q\leq d$, we need to establish some
properties of minimizers of (\ref{constrainedminimization}), and we are
especially interested in their behavior near the boundary. At this point, we assume we are given a family $\{u_\ve\}$
of minimizers for (\ref{constrainedminimization}) and
\begin{equation}
E_\ve(u_\ve)\leq \Lambda:= I_0(d,A)+\pi(|d-p|+|d-q|).
\label{basicbound}
\end{equation}
Also, we suppose also $\ve\leq\ve_0$, where
$\ve_0=\ve_0(\Lambda)>0$ is as in Proposition
\ref{openconstraint}. It follows that maps $u_\ve$ are local
minimizers of $E_\ve(u)$ in $\C{J}$ and therefore they satisfy
(\ref{equation}), (\ref{boundary}).

We will use the notations:  $B_\ve(y)=\{x\in \mathbb{R}^2: |x-y| <\ve\}$,  $\rho_\ve(x)=|u_\ve(x)|$,
$h_\ve(x)$ is the unique solution of (\ref{system_for_r})
(associated to $u_\ve$), and $\C{L}$ is the contour as in Lemma \ref{sensadeg}.
The contour $\C{L}$
separates the two open subdomains $Q^{\pm}$ in $A$, where
$Q^{+}$ is the domain enclosed by $\partial \Omega$ and
$\C{L}$ and  $Q^{-}=A\setminus(Q^{+}\cup \C{L})$.
We also set $Q_\ve^{\pm}=\{x\in Q^{\pm};\,\rho_\ve^2(x)\leq 1-\ve^{1/2}\}$.

\subsection{Proof of Theorem \ref{asymptoticth}}
Since $|\nabla h_\ve|\leq |\nabla u_\ve|$ (by Lemma \ref{maximummodulus}),
the family $\{h_\ve\}$ is bounded in $H^1(A)$, and therefore
there is a sequence $\ve_k\to 0$ such that
\begin{equation}
h_{\ve_k}\to h\quad \text{weakly in $H^1(A)$, as}\ k\to\infty.
\label{convergence_of_r}
\end{equation}
In order to identify $h$, we make use of
Lemma \ref{Lem_Miron} to obtain that, up to a subsequence,
maps $u_{\ve_k}$ converge to a $S^1$-valued map $u$ in
$C^1_{\rm loc}(A)$. Since
$\partial_{x_1} u\times \partial_{x_2} u =0$ a.e. in $A$, we have
$\Delta h_{\ve_k}=2\partial_{x_1}
u_{\ve_k}\times\partial_{x_2} u_{\ve_k}\to 0$
in $C^0_{\rm loc}(A)$, thus $h$ is a
harmonic function. Moreover, $h=1$ on $\partial \Omega$ and
$h={\rm Const}$
on $\partial \omega$.
On the other hand,
$$
{\rm abdeg}(u_{\ve_k})=
\frac{1}{2\pi}\int_A\nabla h_{\ve_k}\cdot\nabla V{\rm d}x \to
\frac{1}{2\pi}\int_A\nabla h\cdot\nabla V{\rm d}x=
\frac{1}{2\pi}\int_{\partial \Omega} \frac{\partial h}{\partial \nu}{\rm d}s.
$$
According to the property c) of ${\rm abdeg}(\,\cdot\,)$
(see Proposition \ref{svoistvoc} in Section \ref{section3}),
${\rm abdeg}(u_\ve)\to d$, as $\ve\to 0$,
therefore $h=h_0$ (where $h_0$ is the
unique solution of (\ref{obschayaotsenka})) and the convergence in
(\ref{convergence_of_r})
holds for the whole family $\{h_\ve\}$. Thus, applying again
Lemma \ref{Lem_Miron}, we obtain
\begin{equation}
h_\ve \to h_0 \quad \text{in $C^1_{\rm loc}$(A)}\
\text{and weakly in $H^1(A)$, as}\ \ve\to 0.
\label{convergence_of_r_precise}
\end{equation}
By (\ref{basicbound}) and Lemma \ref{Lem_Miron}, maps $u_\ve$ converge,
up to a subsequence, to
$u\in H^1(A;S^1)$ in $C^1_{\rm loc}(A)$ and weakly in $H^1(A)$.
Moreover, ${\rm abdeg}(u)=d$ and in view of
(\ref{convergence_of_r_precise})
$|\nabla u|=|\nabla h_0|$ a.e. in $A$. It follows that $u$
is a solution of the minimization problem (\ref{S^1minimization}).

In order to demonstrate the energy convergence stated in Theorem
\ref{asymptoticth}, we argue as follows: by using two pointwise
equalities
$|\nabla u_\ve|^2=2\partial_{x_1} u_\ve\times
\partial_{x_2} u_\ve
+4|\partial_{\bar z} u_\ve|^2$
and $|\nabla u_\ve|^2=-2
\partial_{x_1} u_\ve
\times\partial_{x_2} u_\ve +4|\partial_{z} u_\ve|^2$ and the
pointwise inequality $|\nabla u_\ve|\geq |\nabla h_\ve|$, we have
\begin{multline}
\frac{1}{2}\int_A|\nabla u_\ve|^2{\rm d}x\geq
-\int_{Q_\ve^{+}}
\partial_{x_1} u_\ve\times
\partial_{x_2} u_\ve{\rm d}x
+2\int_{Q_\ve^{+}}|\partial_z u_\ve|^2{\rm d}x\\
+\int_{Q_\ve^{-}}
\partial_{x_1} u_\ve\times
\partial_{x_2} u_\ve{\rm d}x
+2\int_{Q_\ve^{-}}|\partial_{\bar z} u_\ve|^2{\rm d}x
+\frac{1}{2}\int_{A\setminus(Q_\ve^{+}\cup Q_\ve^{-})}|\nabla h_\ve|^2
{\rm d}x.
\label{dlinnayaotsenka}
\end{multline}
Let us estimate the right hand side of (\ref{dlinnayaotsenka}) from below.
Introducing $\sigma_\ve(x)=\max\{\rho_\ve^2(x),1-\ve^{1/2}\}$,
we have (by (\ref{svoistvo_r1})
(\ref{svoistvo_r2}))
\begin{equation}
{\rm div}(\frac{1}{\sigma_\ve(x)}\nabla h_\ve)=
\frac{2}{1-\ve^{1/2}}\begin{cases}
0\quad \text{in}\quad A\setminus(Q_\ve^{+}\cup Q_\ve^{-});\\
\displaystyle
\partial_{x_1} u_\ve\times\partial_{x_2} u_\ve\ \text{otherwise}.
\end{cases}
\label{auxequation}
\end{equation}
Integrating  (\ref{auxequation}) over $Q^{+}$, we get
for sufficiently small $\ve$,
\begin{multline*}
\frac{2}{1-\ve^{1/2}}
\int_{Q_\ve^{+}}\partial_{x_1} u_\ve\times\partial_{x_2} u_\ve{\rm d}x
=
\int_{\partial\Omega}\frac{\partial h_\ve}{\partial \nu}{\rm d}s-
\int_{\C{L}}\frac{\partial h_\ve}{\partial \nu}\frac{{\rm d}s}{\rho_\ve^2(x)}
\\
=\int_{\partial\Omega}u_\ve\times\frac {\partial u_\ve}{\partial \tau}{\rm d}s
-\int_{\C{L}}\frac{u_\ve}{|u_\ve|}\times
\frac{\partial}{\partial \tau}\frac{u_\ve}{|u_\ve|}{\rm d}s=
2\pi(q-d),
\end{multline*}
where we have used Lemma \ref{Lem_Miron} and Lemma \ref{sensadeg}.
Thus, we have
\begin{equation}
\int_{Q_\ve^{+}}
\partial_{x_1} u_\ve\times\partial_{x_2} u_\ve{\rm d}x
=(1-\ve^{1/2})\pi(q-d).
\label{est1}
\end{equation}
Similarly, integrating (\ref{auxequation}) over $Q^{-}$ we obtain
\begin{equation}
\int_{Q_\ve^{-}}
\partial_{x_1} u_\ve\times\partial_{x_2} u_\ve{\rm d}x
=(1-\ve^{1/2})\pi(d-p).
\label{est2}
\end{equation}
In order to estimate the last term in the right hand side of
(\ref{dlinnayaotsenka}), we write it as
\begin{multline*}
\int_{A\setminus(Q_\ve^{+}\cup Q_\ve^{-})}|\nabla h_\ve|^2{\rm d}x
=
\int_{A\setminus(Q_\ve^{+}\cup Q_\ve^{-})}|\nabla h_\ve-\nabla h_0|^2
{\rm d}x \\
+\int_{A}(2\nabla h_\ve-\nabla h_0)\cdot \nabla h_0 {\rm d}x
-\int_{Q_\ve^{+}\cup Q_\ve^{-}}(2\nabla h_\ve-\nabla h_0)\cdot
\nabla h_0 {\rm d}x,
\end{multline*}
and note that by virtue of (\ref{basicbound}) the measure of
$Q_\ve^{+}\cup Q_\ve^{-}$ vanishes as $\ve\to 0$,
so that
\begin{equation}
\int_{A\setminus(Q_\ve^{+}\cup Q_\ve^{-})}|\nabla h_\ve|^2{\rm d}x=
\int_{A\setminus(Q_\ve^{+}\cup Q_\ve^{-})}|\nabla h_\ve-\nabla h_0|^2
{\rm d}x
+\int_{A}|\nabla h_0|^2 {\rm d}x+o(1).
\label{est3}
\end{equation}

Thus (\ref{est1}-\ref{est3}, \ref{dlinnayaotsenka}, \ref{basicbound})
imply $E_\ve(u_\ve)\to E_0(u)+\pi(|d-p|+|d-q|)$. \hfill $\square$

As a byproduct of the above proof by, (\ref{est1})-(\ref{est3}), (\ref{dlinnayaotsenka}) and (\ref{basicbound}) we get
\begin{equation}
\int_A(|u_\ve|^2-1)^2{\rm d}x=o(\ve^2),
\label{est_for_nonlinear_part}
\end{equation}
\begin{equation}
\label{ve>e^1/2}
\int_{A\setminus(Q^+_\ve\cup Q^-_\ve)}|\nabla h_\ve-\nabla h_0|^2
{\rm d}x=o(1),
\end{equation}
\begin{equation}
\label{partialz}
\int_{Q^+_\ve}|\partial_{z} u_\ve|^2{\rm d}x=o(1)\quad
\int_{Q^-_\ve}|\partial_{\bar z} u_\ve|^2{\rm d}x=o(1).
\end{equation}
%We also showed in the course of the proof that
%\begin{equation}
%\label{integraldr}
%1>r_\ve(x)>r_\ve(\partial\omega)\ \text{on}\ \Gamma\quad
%\text{and}\quad
%\int_{\Gamma}
%\frac{\partial r_\ve}{\partial \nu}\frac{d\sigma}{\rho_\ve^2(x)}=2\pi d,
%\end{equation}
%when $\ve$ is small.

\subsection{Properties of minimizers of (\ref{constrainedminimization})
for small $\ve$}

First, by using (\ref{est_for_nonlinear_part}) and the following
methods of \cite{BBH} we get that $\rho_\ve$ converges to $1$
uniformly on compacts in $A$. Moreover, we have

\begin{lem} For any $\mu>0$ we have
\begin{equation}
\sup\{{\rm dist}(y,\partial A);\, y\in A, \rho^2_\ve(y)<1-\mu\}=o(\ve).
\end{equation}
\label{lem_for_modulus}
\end{lem}
\noindent
{\bf Proof.} Assume by contradiction, that for a sequence $\ve_k\to 0$ and $\gamma>0$
we have $\rho^2_{\ve_k}(y_k)<1-\mu$ and $\dist(y_k,\partial A)\geq\gamma \ve_k$.
Due to (\ref{est_2_miron}),
$
|\nabla |u_{\ve_k}|^2|\leq \alpha/\ve_k$ in $B_{\lambda\ve_k}(y_k)$,
where $0<\lambda<\gamma$ and $\alpha(=\alpha(\lambda))$ is independent of $\ve_k$. It follows that $|u_{\ve_k}(x)|^2<1-\mu+\delta\alpha$ when $x\in B_{\delta\ve_k}(y_k)$ and $\delta<\lambda$. Then
$B_{\delta\ve_k}(y_k)\subset A$ and
$$
\frac{1}{\ve^2_k}\int_{B_{\delta\ve_k}(y_k)}(|u_{\ve_k}|^2-1)^2{\rm d}x\geq \pi(\mu-\delta\alpha)^2\delta^2>0,
$$
as soon as $0<\delta<\min\{\lambda, \mu/(2\alpha)\}$. This contradicts (\ref{est_for_nonlinear_part}).\hfill $\square$

Important properties of $u_\ve$ and $h_\ve$, in a vicinity of the
boudary $\partial A$, are established in

\begin{lem}
%\label{structurelemma1}
For any $0<\mu<1$ and $\kappa<1$ there are
$\hat\ve_1(\mu),\hat\ve_2(\mu,\kappa)>0$ such that
if $\rho^2_\ve(y)\leq 1-\mu$  then
\begin{itemize}
%\item[{\rm (i)}\ \ ] ${\rm dist}(y,\partial A)<\ve/4$;
\item[{\rm (a)}] for $\ve<\hat\ve_1(\mu)$ we have $h_\ve(y)\geq 1+\mu/4$
if ${\rm dist}(y,\partial \Omega)<\ve$
and $h_\ve(y)\leq r_\ve(\partial \omega)-\mu/4$ if
${\rm dist}(y,\partial \omega)<\ve$;
\item[{\rm (b)}] for $\ve<\hat\ve_2(\mu,\kappa)$ we have
\begin{equation}
\label{good_otsenka}
\frac{1}{2}\int_{A\cap B_\ve(y)}|\nabla u_\ve|^2 {\rm d}x\geq \kappa\pi.
\end{equation}
\end{itemize}
\label{energy_lem}
\end{lem}
\noindent {\bf Proof.}(by contradiction) Let us assume that either
(a) or (b) is violated for a sequence $\ve_k\to 0$ and some $y=y_k$
such that $h_{\ve_k}(y_k)\leq 1-\mu$. According to Lemma
\ref{lem_for_modulus}, $y_k\to\partial A$. For the definiteness we
suppose that $y_k\to\partial \Omega$, then (by Lemma
\ref{lem_for_modulus})
\begin{equation}
\label{odist}
{\rm dist}(y_k,\partial \Omega)=o(\ve_k).
\end{equation}
Let $u_{\ve}$ be continued in $\omega$ in such a way that
$\|u_\ve\|_{H^1(\Omega)}\leq
C\|u_\ve\|_{H^1(A)}$ and $|u_\ve|\leq 1$ in $\Omega$,
where $C$ is independed of $\ve$. We also assume that $h_\ve=h_\ve(\partial\omega)$ on
$\omega$.
Following \cite{BM2} we rescale $u_{\ve_k}$ and
$h_{\ve_k}$ by a conformal map that 'moves' $y_k$ away from the
boundary. Fix a conformal mapping $\eta$ from
$\Omega$ onto the unit disk $B_1(0)$.
%Then the sequences
%$u_{\ve_k}(\eta^{-1}(\ \cdot\ ))$,  $r_{\ve_k}(\eta^{-1}(\ \cdot\ ))$
%are bounded in $H^1(B_1(0))$.
We introduce the conformal map
$\zeta_k(z)=(z-\eta(y_k))/(\bar\eta(y_k)z-1)$ from $B_1(0)$ onto
itself and set $U_k(z)=u_{\ve_k}(\eta^{-1}(\zeta_k(z)))$,
$H_k(z)=h_{\ve_k}(\eta^{-1}(\zeta_k(z)))$. It is easy to see that
$\|U_k\|_{H^1(\Omega)}\leq C$ and $\|H_k\|_{H^1(\Omega)}\leq C$
with some $C$ independent of $k$. Therefore, without loss of
generality, we can assume that $U_k$ and $H_k$ $H^1$-weakly
converge to limits $U$ and $H$, respectively, as $k\to\infty$.

Arguing, as in \cite{BM2}(Section 4), we can show that $U_k\to U$ in
$C^1_{\rm loc}(B_1(0))$ and that $\Delta U=0$ in $B_1(0)$.
Therefore, $|U(0)|^2=\lim_{k\to\infty} |U_k(0)|^2=
\lim_{k\to\infty} |u_{\ve_k}(y_k)|^2\leq 1-\mu$.
We also have $|U|=1$ a.e. on $S^1$. Show now that
$\partial_z U=0$ in $B_1(0)$. Indeed, by the maximum principle $|U|<1$ in $B_1(0)$,
hence $\max_{B_t(0)}|U_k(z)|^2<1-\ve_k^{1/2}$ for any fixed $0<t<1$ and sufficiently large $k$. It follows that for such $k$
we have $\eta^{-1}(\zeta_k(B_t(0)))\subset Q^+_\ve$. Then in view of (\ref{partialz}) we get, by using
the conformality of the maps  $\eta^{-1}$ and $\zeta_k$,
$$
\int_{B_t(0)}|\partial_z U_k|^2dx=\int_{\eta^{-1}(\zeta_k(B_t(0)))}|\partial_z u_{\ve_k}|^2dx
\leq\int_{A^+_\ve}|\partial_z u_{\ve_k}|^2dx\to 0
$$
This implies that  $\partial_z U=0$ in $B_1(0)$.

In order to show (b) we use the pointwise equalities
$\frac{1}{2}|\nabla U|^2=-\partial_{x_1} U\times \partial_{x_2}
U+\frac{1}{4}|\partial_z U|^2$ and $\partial_z U=0$ to obtain
$$
\frac{1}{2}\int_{B_1(0)}|\nabla U|^2{\rm d}x=
-\int_{B_1(0)}\partial_{x_1} U\times \partial_{x_2}U{\rm d}x=
-\pi{\rm \deg}(U,S^1).
$$
As $U\not\equiv {\rm Const}$, we therefore have
$\frac{1}{2}\int_{B_1(0)}|\nabla U|^2{\rm d}x\geq \pi$.
It follows that there is $0<t<1$ such that
\begin{equation}
\label{zvezda}
\frac{1}{2}\int_{B_t(0)}|\nabla U|^2dx>\kappa \pi.
\end{equation}
The image $\zeta_k(B_t(0))$ of the disk $B_t(0)$ is the disk
$B_{t_k}(\xi_k)$ with the radius $t_k=\frac{t(1-|\eta(y_k)|^2)}{1-t^2|\eta(y_k)|^2}$
and the center at $\xi_k=\frac{1-t^2}{1-t^2|\eta(y_k)|^2}$. According to
(\ref{odist}) $t_k=o(\ve_k)$ for $k\to\infty$, hence
$\eta^{-1}(B_{t_k}(\xi_k))\subset B_{\ve_k}(y_k)$ when $k$ is sufficiently large.
Then, by using the conformal invariance
and lower semicontinuity of the Dirichlet integral,
and bound (\ref{zvezda}), we get
$$
\int_{B_{\ve_k}(y_k)}|\nabla u_{\ve_k}|^2{\rm d}x\geq
\int_{\eta^{-1}(\zeta_k(B_t(0)))}|\nabla u_{\ve_k}|^2{\rm d}x=
\int_{B_t(0)}|\nabla U_k|^2{\rm d}x>2\kappa \pi \ \text{as}\ k\to\infty.
$$

In order to show that $h_{\ve_k}(y_k)=H_k(0)>1+\mu/4$ when $k\to\infty$ we note
that the system (\ref{system_for_r}) is conformally invariant, i.e.
$$
\nabla^\bot H_k=( U_k\times\partial_{x_1} U_k, U_k\times\partial_{x_2}U_k)
\quad\text{in}\quad \zeta_k^{-1}(\eta(A)).
$$
Then, bearing in mind the convergence properties of $U_k$, we obtain that $H_k\to H$ in $C^1_{\rm loc}(B_1(0))$
and
$$
\nabla^\bot H=(U\times \partial_{x_1} U , U\times\partial_{x_2}U)
=-\frac{1}{2}\nabla^\bot(|U|^2)
\quad\text{in}\quad B_1(0),
$$
where we have used the fact that $\partial_z U=0$ in $B_1(0)$.
Since $H=|U|=1$ on $\partial B_1(0)$ we have $H=\frac{3}{2}-\frac{1}{2}|U|^2$ in $B_1(0)$,
therefore
$$
\lim_{k\to \infty} h_{\ve_k}(y_k)=\lim_{k\to \infty} H_k(0)=\frac{3}{2}-\frac{1}{2}|U(0)|^2
\geq 1+\frac{\mu}{2}.
$$
Lemma is proved.\hfill $\square$

\begin{rem}
\label{remvortexless}
Lemma \ref{energy_lem} implies that in the case
when $p=q=d$ minimizers of (\ref{constrainedminimization}) are
vortexless for sufficiently small $\ve$. Really, by Theorem
\ref{asymptoticth} they
%minimizers $\{u_\ve\}$
$H^1$-strongly converge, up to a subsequence, as $\ve\to 0$ to a minimizing harmonic map $u \in \C{J}^{(d)}_{dd}$. On the other
hand  (\ref{good_otsenka}), exhibits the energy concentration
property near zeros of minimizers, which is incompatible with
the strong $H^1$-convergence.
\end{rem}

The following result, describing the structure of the function
$h_\ve$ for small $\ve$ plays a crucial role in the proof of the
main technical result (Lemma \ref{mainotsenka}) in Section
\ref{section7}.

\begin{lem}
\label{sructure}
We have, for small $\ve$, $\ve<\ve_4$ (where $\ve_4>0$),
\begin{itemize}
\item[(i)\ ] $\rho^2_\ve(x)\geq 1/2$ when
$h_\ve(\partial \omega)-1/8\leq h_\ve(x)\leq 9/8$
and $h_\ve(\partial \omega)<\min_{\C{L}}h_\ve(x)\leq\max_{\C{L}}h_\ve(x)<
h_\ve(\partial\Omega)$, while
if $\rho^2_\ve(x)<1/2$ then either
$$
{\rm dist}(x,\partial\Omega)<{\rm dist}(\C{L},\partial\Omega)\ \text{and}\
h_\ve(x)>9/8$$
or
$${\rm dist}(x,\partial\omega)<{\rm dist}(\C{L},\partial\omega)\ \text{and}\
h_\ve(x)<h_\ve(\partial \omega)-1/8;$$
%\item[\ (ii)] for any regular value  $t\in (r_\ve(\partial \omega), 1)$
%of $r_\ve$
%the level set $r_\ve(x)=t$ is a smooth connected curve
%(without selfintersecting) that encloses $\omega$,
%this curve being oriented counterclockwise,
%\begin{equation}
%\label{integraluu}
%\frac{1}{2\pi}\int_{r_\ve(x)=t}
%u_\ve\times \frac{\partial u_\ve}{\partial\tau}d\sigma=d,
%\end{equation}
%and $u_\ve\times \frac{\partial u_\ve}{\partial\tau}>0$;
\item[(ii)] there are
$x^*_\ve\in \partial\Omega$, $x^{**}_\ve\in \partial\omega$ such that
$\frac{\partial h_\ve}{\partial \nu}(x^*_\ve)>0$ and
$\frac{\partial h_\ve}{\partial \nu}(x^{**}_\ve)>0$.
\end{itemize}
\end{lem}
\noindent
{\bf Proof.} (i) follows from Lemma \ref{lem_for_modulus}, Lemma \ref{energy_lem}
and the convergence properties of $h_\ve$ as $\ve\to 0$ established in the course of
the proving Theorem \ref{asymptoticth}.
%By Sard's lemma a.e. $t\in(r_\ve(\partial\omega),1)$
%is a regular value of $r_\ve$. Then, according to (\ref{integraldr}), we can chose
%a regular value $t^*$ satisfying $1>t^*>\max_{x\in\Gamma} r_\ve(x)$. The level set
%$\{x:r_\ve(x)=t^*\}$ is a one dimensional $C^1$ manifold. It consists of a finite number
%of connected components that all are contained in $A$. Assume that there is a component
%not enclosing $\omega$. Chose a minimal one $\Gamma^*$ in the sense that it encloses no other components. Then we have either  domain enclosed by $\Gamma^*$  is
To prove (ii) we argue as follows.
Let $\hat\ve_2(\mu,\kappa)$ be the best (biggest) constant
in Lemma \ref{energy_lem}. Then $\hat\ve_2(\mu,\kappa)$ is
increasing in $\mu$ and decreasing in $\kappa$. For $k=1,2,\dots$ set
$$
\hat \mu_\ve={1}/{k}\quad \text{when}\quad
\min\{\hat\ve_2({\textstyle \frac{1}{k+1}},{\textstyle\frac{k}{k+1}}), {1}/({k+1})\}
\leq \ve
< \min\{\hat\ve_2({\textstyle \frac{1}{k}},{\textstyle \frac{k-1}{k}}), {1}/{k}\}.
$$
Then $\hat\mu_\ve\to 0$ as $\ve\to 0$ and (\ref{good_otsenka})
is satisfied with $\kappa=1-\hat\mu_\ve$ when $\rho^2_\ve(y)<1-\hat\mu_\ve$;
the same being true when $\hat\mu_\ve$ is
replaced by $\mu_\ve=\max\{\hat\mu_\ve,\ve^{1/2}\}$.
We pick a point $x_\ve^{(1)}$ in $A$ such that $\rho^2_\ve(x_\ve^{(1)})<1-\mu_\ve$; then we pick a point $x_\ve^{(2)}$  in
$A\setminus B_{2\ve}(x_\ve^{(1)})$ such that
$\rho^2_\ve(x_\ve^{(1)})<\mu_\ve$ , e.t.c. unless for some $K_\ve$
we have $\rho^2_\ve(x)\geq 1-\mu_\ve$ on
$A\setminus \cup_{k=1}^{K_\ve} B_{2\ve}(x_\ve^{(k)})$.
By the construction of $\mu_\ve$, since disks $B_{\ve}(x_\ve^{(k)})$ are disjoint,
$$
\frac{1}{2}\int_A|\nabla u^\ve|^2dx\geq \frac{1}{2}\sum_1^{K_\ve}\int_{A\cap B_{\ve}(x_\ve^{(k)})}|\nabla u^\ve|^2dx
\geq  K_\ve(1-\mu_\ve)\pi.
$$
Therefore by (\ref{basicbound}) we have a uniform bound $K_\ve\leq C$.
Arguing as in \cite{BBH} (Chapter IV, Theorem IV.1) we can
increase the radii of disks
to $\ve\lambda>2\ve$ (with $\lambda$ independent of $\ve$) and take a subset
$I_\ve$ of $\{1,\dots,K_\ve\}$ such that
$$
\cup_{k\in I_\ve} B_{\ve\lambda}(x_\ve^{k})\supset\cup_{k=1}^{K_\ve} B_{2\ve}(x_\ve^{k})
\quad\text{and}\quad \dist(x_\ve^{k^\prime},x_\ve^{k})>4\ve\lambda
\quad\text{for different}\quad k,k^\prime\in I_\ve.
$$
We also have
$$
\rho^2_\ve(x)\geq 1-\mu_\ve\ \text{on}\ A\setminus
\cup_{k\in I_\ve} B_{\ve\lambda}(x_\ve^{k})
$$

Assume $h_\ve$ and $h_0$
extend to the whole $\D{R}^2$ and set
$h_\ve=h_\ve(\partial \omega)$, $h_0=h_0(\partial \omega)$  in $\omega$,
and $h_\ve=h_0=1$ in $\D{R}^2\setminus \Omega$. Since
$\rho^2_\ve(x)\geq 1-\mu_\ve \geq 1-\ve^{1/2}$ in
$D_\ve^{(k)}
=B_{2\lambda\ve}(x_\ve^{(k)})\setminus B_{\lambda\ve}(x_\ve^{(k)})$,
by (\ref{ve>e^1/2}) we have
$$
\int_{D_\ve^{(k)}}
|\nabla h_\ve|^2{\rm d}x\leq
2\int_{D_\ve^{(k)}}
(|\nabla (h_\ve-h_0)|^2+|\nabla h_0|^2){\rm d}x
\to 0\ \text{as}\ \ve\to 0.
$$
Then, writing the integral over $D_\ve^{(k)}$ in the polar coordinates
with the center at $x_\ve^{(k)}$ and using
Fubini's theorem, we can find $\lambda_\ve^{(k)}$,
$\lambda\ve\leq\lambda_\ve^{(k)}\leq 2\lambda\ve$, such that
$$
\int_{|x-x_\ve^{(k)}|=\lambda_\ve^{(k)}}
|\nabla h_\ve|^2d\sigma=o(1/\ve).
$$
Therefore, by the Cauchy-Shwartz inequality,
\begin{equation}
\int_{|x-x_\ve^{(k)}|=\lambda_\ve^{(k)}}
|\nabla h_\ve|{\rm d}s\leq (\pi\lambda_\ve^{(k)})^{1/2}\left\{
\int_{|x-x_\ve^{(k)}|=\lambda_\ve^{(k)}}
|\nabla h_\ve|^2{\rm d}s\right\}^{1/2}=o(1).
\label{ogo}
\end{equation}
Now, integrate(\ref{svoistvo_r2}) over
$Q^{+}\setminus\cup_{k\in I_\ve} B_{\lambda_\ve^{(k)}}(x_\ve^{(k)})$ to get,
according to (\ref{ogo}) and Lemma \ref{sensadeg},
$$
\int_{\Gamma_\ve}\frac{\partial h_\ve}{\partial \nu}{\rm d}s=2\pi d
+\sum_{k\in I_\ve}
\int_{|x-x_\ve^{(k)}|=\lambda_\ve^{(k)}}
\frac{\partial h_\ve}{\partial \nu}\frac{{\rm d}s}{\rho_\ve^2}=2\pi d - o(1)\ \text{when}\ \ve\to 0,
$$
where $I^\prime_\ve$ denotes the subset of indexes $k\in I_\ve$ such that
$B_{\lambda_\ve^{(k)}}(x_\ve^{(k)})\cap Q^{+}\not=\emptyset$, and
$\Gamma_\ve=\partial\Omega\setminus\cup_{k\in I^\prime_\ve} B_{\lambda_\ve^{(k)}}(x_\ve^{(k)})$.
Therefore, there is $x_\ve^{*}\in \Gamma_\ve$ such that $\frac{\partial h_\ve}{\partial \nu}(x_\ve^{*})>0$. Similarly, we can show that
on $\gamma_\ve=\partial\omega\setminus\cup_{k\in I_\ve\setminus I_\ve^\prime} B_{\lambda_\ve^{(k)}}(x_\ve^{(k)})$ there is $x_\ve^{**}$
such that $\frac{\partial h_\ve}{\partial \nu}(x_\ve^{**})>0$.
\hfill $\square$

\section{Inductive proof of Theorem \ref{mainth}}

\label{section7}
Fix $\Lambda> 0$ and an integer $d>0$ such that
$\Lambda>I_0(d,A)$,
and let $\ae_0$ be the greatest integer such that
$$
I_0(d,A)+\pi\ae_0<\Lambda.
$$
Clearly, $\ae_0\geq 0$. In this Section we show that for small $\ve$ the
infimum in problem (\ref{constrainedminimization}) is always attained,
provided integers $p$, $q$ satisfy
\begin{equation}
\label{boundforpq}
%\Lambda>\frac{2(\pi d)^2}{{\rm cap}(A)}+\pi(|d-p|+|d-q|) \quad\text{and}\quad
%p\leq d,\ q\leq d.
\ae(p,q)\leq\ae_0 \quad\text{and}\quad
p\leq d,\ q\leq d,
\end{equation}
where $\ae(p,q)=|d-q|+|d-p|$.

\begin{prop}
%{\rm (induction hypothesis)}
\label{mainstatement} Given an integer $K\leq\ae_0$. Let $p\leq q$,
$q\leq d$ be integers such that $\ae(p,q)\leq K$.
%satisfying (\ref{boundforpq}).
Then,
for sufficiently small $\ve$ the infimum in problem (\ref{constrainedminimization})
is always attained and
\begin{equation}
\label{induction}
m_\ve(p,q,d)\leq m_\ve(l,m,d)+\pi((l-p)+(m-q))\quad \text{when}\quad p\leq l\leq d,\ q\leq m\leq d.
\end{equation}
Moreover, the inequality in (\ref{induction}) is strict unless $l=p$ and $m=q$.
\end{prop}
\noindent {\bf Proof.} The proof is by induction on $K$. The basis
of induction (K=0) is established in Section \ref{section4} (cf.
Lemma \ref{firststep}). The demonstration of the induction step
relies on the following Lemma, whose proof is in the end of this
section.
\begin{lem}
\label{mainotsenka}
Assume integers $p$ and $q$ satisfy (\ref{boundforpq}), and
for $\ve<\ve_5$, $\ve_5>0$, there exists a minimizer
$u_\ve$ of problem (\ref{constrainedminimization}) whose GL energy
$E_\ve(u_\ve)$ satisfies the bound (\ref{basicbound}).
%Let $u_\ve$ be a minimizer of $E_\ve(u)$ in ${\cal J}^{(d)}_{pq}$ (integers $p$, $q$ satisfy
%assumptions of Proposition \ref{mainstatement})
%and let $\ve$ be so small
Then for any $\ve<\min\{\ve_4,\ve_5\}$ (where $\ve_4$ is as in Lemma \ref{sructure}) there exists
$v_\ve\in {\C{J}}_{p (q-1)}$ such that
%{\rm (a)}
${\rm abdeg}(v_\ve)\in(d-1/2,d+1/2)$ and
%{\rm (b)}
\begin{equation}
\label{samayaglavnayaotsenka} E_\ve(v_\ve)<m_\ve(p,q,d)+\pi.
\end{equation}
Similarly, in ${\C{J}}_{(p-1)q}$ there exists a testing map
(still denoted $v_\ve$) satisfying (\ref{samayaglavnayaotsenka}),
and such that ${\rm abdeg}(v_\ve)\in(d-1/2,d+1/2)$.
\end{lem}

In view of Lemma \ref{mainotsenka}, to prove the claim of
Proposition \ref{mainstatement} for $K+1$ in place of $K$ it
suffices to show that $m_\ve(p,q-1,d)$ is always attained for
sufficiently small $\ve$ whenever $p\leq d$, $q\leq d$ and
$\ae(p,q)=K$ (the attainability of $m_\ve(p-1,q,d)$ is proved
similarly). Let $u$ be a weak $H^1$-limit of a minimizing sequence
$(u^{(k)})$. According to Lemma \ref{mainotsenka}, such a
minimizing sequence exists, moreover, this lemma, with the induction
hypothesis, imply
\begin{equation}
\label{promezhutochnoe}
\limsup_{k\to\infty} E_\ve(u^{(k)})<
m_\ve(l^\prime,m^\prime,d)+\pi(l^\prime-p)+\pi(m^\prime-(q-1)),
\end{equation}
where $p\leq l^\prime\leq d$, $q-1\leq m^\prime\leq d$ and
$\ae(l^\prime,m^\prime)\leq K$.  We know that ${\rm
abdeg}(u)=\lim_{k\to\infty}{\rm abdeg}(u^{(k)})\in (d-1/2,d+1/2)$
when $\ve<\ve_0$ (where $\ve_0=\ve_0(\Lambda)$ is as in
Proposition \ref{openconstraint}), hence $u$ is a solution of the
GL equation (\ref{equation}) (see arguments in Lemma
\ref{firststep}). Therefore, if we write $\liminf_{k\to\infty}
E_\ve(u^{(k)})\leq \limsup_{k\to\infty} E_\ve(u^{(k)})$ and apply
successively Lemma \ref{pricelemma} and Lemma
\ref{roughlowerbound} to the left hand side, we get by using
(\ref{promezhutochnoe}) with $l^\prime=m^\prime=d$ that for $\ve$
sufficiently small
%\begin{multline}
$$
%\label{posledsravnen}
I_0(d,A)-\frac{\pi}{2}+\pi\ae(l,m)+\pi(|l-p|+|m-(q-1)|)
\leq m_\ve(d,d,d)
+\pi\ae(p,q-1),
$$
%\end{multline}
where $l={\rm deg}(u,\partial\omega)$, $m={\rm deg}(u,\partial\Omega)$.
Thanks
to Lemma \ref{ochenprostayalemma} $m_\ve(d,d,d)\leq I_0(d,A)$, thus
$$
|l-d|+|l-p|\leq |p-d|+1/2\ \text{and}\
|m-d|+|m-(q-1)|\leq |(q-1)-d|+1/2.
$$
Since $l$ and $m$ are integers, it follows that $p\leq l\leq d$,
$q-1\leq m\leq d$. Now, assuming $l\not=p$ or $m\not=q-1$, we
use Lemma \ref{pricelemma} and (\ref{promezhutochnoe}) with
$l^\prime=l$, $m^\prime=m$ to obtain the following
$$
E_\ve(u)+\pi(l-p)
\leq\liminf_{k\to\infty} E_\ve(u^{(k)})
<m_\ve(l,m,d)
+\pi(l-p)+\pi(m-(q-1)).
$$
On the other hand, $u\in\C{J}_{lm}$ and ${\rm
abdeg}(u)\in[d-1/2,d+1/2]$. Hence $E_\ve(u)\geq m_\ve(l,m,d)$, which
is a contradiction. Therefore $l=p$ and $m=q-1$, i.e. $u$ is an
admissible testing map in problem (\ref{constrainedminimization})
and thus the infimum $m_\ve(p,q-1,d)$ is always attained.\hfill
$\square$

\noindent
{\bf Proof of Lemma \ref{mainotsenka}.}
For simplicity we drop subscript $\ve$. The underlying idea
is to modify the minimizer $u$ of (\ref{constrainedminimization})
in a neighborhood of $\partial \Omega$ as in Proposition \ref{testingmap} (see Section \ref{section5}). In general $u$ is with zeros
now, thus the arguments need to be more sophisticated. Loosely speaking, we construct a testing
map $v$ with an additional "vortex" located "near" $x^*$, where $x^*$
is a point on $\partial\Omega$ such that
$\frac{\partial h}{\partial \nu}(x^*)>0$ (cf. Lemma \ref{sructure}).

{\bf Step 1: Domain decomposition.} Let $1-\delta$, where $\delta>0$, be a
regular value of $h$ (thanks to Stard's lemma
this holds for almost all $\delta$). Consider the subdomain
of $A$ where $h>1-\delta$.
There is a (unique) connected component $D_\delta$ of this subdomain,
such that $\partial D_\delta\supset\partial\Omega$. Since
$h(\partial\Omega)=1>h(\partial\omega)$, when $\delta$
is sufficiently small the boundary of
$D_\delta$ contains a connected component $\Gamma_\delta\not=\partial\Omega$
enclosing $\omega$. According to Lemma \ref{sructure}, we can choose $\delta$ small enough $\delta<\delta_0$ ($\delta_0>0$)
so that the domain enclosed by $\Gamma_\delta$ and
$\partial\Omega$ lies away from the contour $\C{L}$, i.e.
$\Gamma_\delta$ also encloses $\C{L}$, moreover if $\rho(x)(=|u(x)|)$
vanishes at a point $x$ of this domain, then $h(x)>1$.
Therefore, the minimum of $h$ over the closure of the
forementioned domain cannot be attained at any interior point,
otherwise $h$ satisfies ${\rm div}(\frac{1}{\rho^2} h)=0$
in a neighborhood of this point, which is impossible. In other words,
$h>1-\delta$ in the domain enclosed by $\Gamma_\delta$ and
$\partial\Omega$, i.e. this domain coincides with $D_\delta$.
Thus, the boundary of $D_\delta$ consists of exactly two connected components
$\partial\Omega$ and $\Gamma_\delta$.
%The boundaries of the domains
%$D_{\delta^\prime}$ have the same structure when $0<\delta^\prime<\delta$
%and $1-\delta^\prime$ is a regular value of $h$.
Also, possibly
choosing smaller $\delta_0$, we have that the set
$P=\{x\in D_{\delta};\, h(x)\geq 1\}$ is independent of
$\delta$ (recall that $\delta<\delta_0$ and $1-\delta$ is a regular value of $h$).
Indeed, consider the
set $S_{\delta}=\{x\in A; h(x)>\alpha\}\cap D_{\delta}$
($0<\delta<\delta_0$), where $\alpha$ is a regular value of $h$
and $1<\alpha<9/8$. $S_{\delta}$ consists of a finite
number $n(\delta)$ of connected components. Since
$D_{\delta}\supset D_{\delta^{\prime}}$ if
$\delta>\delta^{\prime}$, the function
$n(\delta)$ is nondecreasing, hence
$n(\delta)=\lim_{\delta^{\prime}\to 0}n(\delta^{\prime})$
when $0<\delta<\delta_0$ (for some $\delta_0>0$) and
$S_{\delta}=S_{\delta^{\prime}}$ if
$\delta,\delta^{\prime}\in(0,\delta_0)$. It follows
that $1-\delta<h<\alpha$ in $D_{\delta}\setminus D_{\delta^{\prime}}$
when $0<\delta^{\prime}<\delta<\delta_0$. For such
$\delta$ and $\delta^{\prime}$ the function $h$ satisfies
${\rm div}(\frac{1}{\rho^2} h)=0$ in $D_{\delta}\setminus D_{\delta^{\prime}}$
(by Lemma \ref{sructure}) while $h<1$ on the boundary of
$D_{\delta}\setminus D_{\delta^{\prime}}$, hence $h<1$
in $D_{\delta}\setminus D_{\delta^{\prime}}$. We see now that
$\{x\in D_{\delta}; h(x)\geq 1\}=P:=\cap_{\delta^\prime<\delta_0} D_{\delta^\prime}$,
when $0<\delta<\delta_0$, as required.

Thus we have, for $\delta<\delta_0$
%Now, according to the structure of $D_\delta$ and $P$ in conjunction
%with the fact that
$$
d^\prime:=
\frac{1}{2\pi}\int_{\Gamma_\delta}\frac{\partial h}{\partial \nu}
\frac{{\rm d} s}{\rho^2}=
\frac{1}{2\pi}\int_{\Gamma_\delta}
\frac{u}{|u|}\times\frac{\partial}{\partial\tau}\frac{u}{|u|}{\rm d}s
={\rm deg}(u,\Gamma_\delta)>0
$$
($1-\delta$ is a regular value of $h$ and $h>1-\delta$  in  $D_\delta$)
and the integer $d^\prime$ is independent of $\delta$. Therefore $u$ admits the representation $u=\rho e^{id^\prime\theta}$
in $G_\delta=D_\delta\setminus P$, where $\theta:G_\delta\to\D{R}\setminus2\pi\D{Z}$
is a smooth function and $\rho>0$ in $\overline G_\delta$.

Without loss of generality we can
assume that $u(x^*)=1$. Since $\nabla h =-d^{\prime}\rho^2\nabla^\bot \theta$
in $G_\delta$ and $\frac{\partial h}{\partial \nu}(x^*)>0$,
the map $x\to(h,\theta)$ from a neighborhood
$G_\delta^\prime$ of $x^*$ onto its image $h(G_\delta^\prime)$
is a $C^1$-diffemorphism. Choosing $\delta$  small enough, we can
assume that $G_\delta^\prime$ is defined by
$$
x\in G_\delta^\prime\ \Longleftrightarrow\ x\in G_\delta,
\ 1-\delta<h(x)<1,\
\theta(x)\in (-\delta,\delta)({\rm mod}\, 2\pi\D{Z}).
$$

\begin{figure}[!h]
\begin{center}
    \psfig{figure=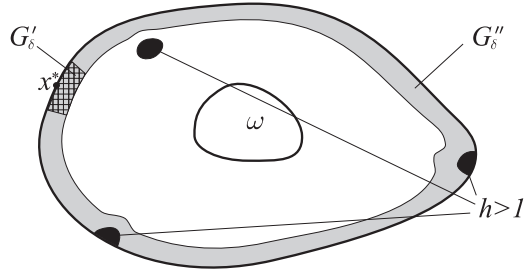,height=1.4in,silent=}
%%%%%%%%%%change bf\psfig{figure=DomDec1.eps,height=2.0in,silent=}
\end{center}
\caption{Domain decomposition.}\label{fig:notations}
\end{figure}

Now we have
$A=G_\delta^\prime\cup G_\delta^{\prime\prime}\cup(A\setminus G_\delta)$
(see Fig. \ref{fig:notations}), where $G_\delta^{\prime\prime}=G_\delta\setminus G_\delta^\prime$.

%Note that if $v\in H^1$ then $|w_t|=1$ on $\partial G_\delta$.

{\bf Step 2: Construction of the testing map.} We seek the tasting map $v$ in
the form
\begin{equation}
\label{test}
v=
\begin{cases}
u\quad \text{in}\ A\setminus G_\delta,\\
\rho w_t\quad \text{in}\ G_\delta,
\end{cases}
\end{equation}
with (unknown for the moment) $w_t=w_t(h(x),\theta(x))$.
Impose the following boundary conditions on $\partial G_\delta$,
\begin{equation}
\label{boundary1}
w_t=
{\rm e}^{id^\prime\theta}
\frac{{\rm e}^{-i\theta}-(1-t\varphi(\theta))}
{{\rm e}^{-i\theta}(1-t\varphi(\theta))-1}
\quad \text{on}\quad \partial G_\delta\setminus\Gamma_\delta
\end{equation}
\begin{equation}
\label{boundary2}
w_t={\rm e}^{id^\prime\theta}\quad \text{on}
\quad \Gamma_\delta,
\end{equation}
where $0\leq\varphi\leq 1$ is a smooth $2\pi$-periodic cut-off
function such that $\varphi(\theta)=1$,
when $\theta\in(-\delta/2,\delta/2)$ (${\rm mod} \ 2\pi\D{Z}$) and $\varphi(\theta)=0$ if $\theta\not\in (-\delta,\delta)+2\pi\D{Z}$.
It is easy to see
%from (\ref{test})
that if $w_t$ (considered as a function of $h,\theta$)
satisfies (\ref{boundary1}) and (\ref{boundary2}), when $h=1$ and $h=1-\delta$, respectively, and is a smooth $2\pi-$periodic
in $\theta$ map defined in the strip $1-\delta\leq h\leq 1$ then
(\ref{test}) defines for $0<t<1$
a map $v\in H^1(A;\D{R}^2)$ such that $|v|=1$ on
$\partial A$ and
\begin{equation}
\deg(v,\partial\Omega)=q-1,\quad \deg(v,\partial\omega)=p.
\label{stepeni}
\end{equation}

Expand the right hand side of (\ref{boundary1}) into the series
\begin{equation}
\label{Fourircoef}
{\rm e}^{id^\prime\theta}
\frac{{\rm e}^{-i\theta}-(1-t\varphi(\theta))}
{{\rm e}^{-i\theta}(1-t\varphi(\theta))-1}
=(1-tb_{-1}(t)){\rm e}^{id^\prime\theta}+
t\sum_{k\not=-1} b_k(t){\rm e}^{-i(k-d^\prime+1)\theta},
\end{equation}
and set
\begin{equation}
\label{formofw}
w_t(h,\theta)=(1-tb_{-1}(t)f_{-1}(h)){\rm e}^{id^\prime\theta} +
t\sum_{k\not=-1} b_k(t) f_{k}(h){\rm e}^{-i(k-d^\prime+1)\theta},
\end{equation}
where functions $f_k$ are defined by (\ref{deffk}) with
$\displaystyle
k_{\pm}=\pm\frac{1}{d^\prime}\sqrt{(k-d^\prime+1)^2+
\lambda-{d^\prime}^2}$.
The positive parameters $t<1$ and $\lambda \geq 2{d^\prime}^2$
are to be specified later on.
%\in C^1([1-\delta,1];\D{C})$ satisfying
%\begin{equation}
%\label{conditionforfk}
%f_k(1-\delta)=0,\ f_k(1)=1.
%\end{equation}
In what concerns the coefficients $b_k(t)$,
%by the standard facts of
%the Fourier analysis
we have
\begin{equation}
\label{asymtcoef}
|b_k(t)-c_k(t)|\leq C(1+|k|)^{-n}, \ \forall n>0,
\end{equation}
where $C=C(n)$ is independent of $t$, and
%$t c_{-k-1}(t)$ for $k\not=-1$ and
%$1-tc_{-1}$ for $k=0$ is the $k$-th Fourier coefficient of
%$({\rm e}^{-i\theta}-(1-t))/
%({\rm e}^{-i\theta}(1-t)-1)$, i.e.
$c_k=(t-2)(1-t)^k$ for $k\geq 0$, $c_{-1}=1$, $c_k=0$ for $k<-1$.
The estimate (\ref{asymtcoef}) is obtained in a standard way, by
comparing the Fourier coefficients in (\ref{Fourircoef}) with
those of ${\rm e}^{id^\prime\theta} ({\rm e}^{-i\theta}-(1-t))/
({\rm e}^{-i\theta}(1-t)-1)$.

% can be written
%as $b_k=(t-2)(1-t)^k+c_k(t)$ for $k\geq 0$, $b_{-1}=1+c_{-1}(t)$ and
%$b_k=c_k(t)$ for $k\leq -2$, where $|c_k(t)|<C(1+|k|)^{-n}$
%with a constant $C=C(n)$ independent of $t$.

%$$
%b_k(t)=
%\begin{cases}
%,\quad \text{if}\quad k=-1\\
%(t-2)(1-t)^k+c_k(t),\quad \text{if}\quad k\geq 0\\
%c_k(t),\quad \text{if}\quad k\leq -2,
%\end{cases}
%$$
%where
%\begin{equation}
%\label{coefficienty}
%c_k(t)=\frac{1}{2\pi}\int_{-\pi}^{\pi}\frac{{\rm e}^{id\theta}}{t}
%\left(
%\frac{{\rm e}^{-i\theta}-(1-t\varphi(\theta))}
%{{\rm e}^{-i\theta}(1-t\varphi(\theta))-1}
%-
%\frac{{\rm e}^{-i\theta}-(1-t)}
%{{\rm e}^{-i\theta}(1-t)-1}
%\right){\rm e}^{i(k-d+1)\theta}d\theta.
%\end{equation}

{\bf Step 3: Verification of (\ref{samayaglavnayaotsenka}).}
Thanks to Lemma \ref{comparisonlemma} we have (since $|w_t|=1$ on
$\partial G_\delta$ due to (\ref{boundary1}, \ref{boundary2}))
$$
E_\ve(v)=E_\ve(u)+L_\ve^{(d^\prime)}(w_t,G_\delta)
$$
where the functional $L_\ve^{(d^\prime)}(w)$ is defined as in (\ref{Lepsilon}).
Let us show, that for sufficiently small $t$
\begin{equation}
\label{dveotsenki1}
L_\ve^{(d^\prime)}(w_t,G_\delta)<\pi.
\end{equation}
To this end, like in the proof of Proposition \ref{testingmap}, consider the quadratic functional
\begin{equation}
\label{functMprime}
M^\prime_\lambda(w_t)=\frac{1}{2}\int_{G_\delta}
({d^\prime}^2|\partial_r w_t|^2+|\partial_\theta w_t|^2 +
\lambda|w_t-{\rm e}^{i\theta}|^2-{d^\prime}^2|w_t|^2)\rho^2 |\nabla \theta|^2{\rm d}x.
\end{equation}
(One can actually show that $w_t$ minimizes functional (\ref{functMprime})
%$M^\prime_\lambda(w)$
under the boundary conditions (\ref{boundary1}), (\ref{boundary2}).) We have, since
$\nabla h =-d^{\prime}\rho^2\nabla^\bot \theta$ in $G_\delta$ and $\rho\leq 1$
in $A$,
\begin{equation}
\int_{G_\delta}\rho^2|\nabla w_t|^2{\rm d}x\leq
\int_{G_\delta}({d^\prime}^2|\partial_r w_t|^2+|\partial_\theta w_t|^2)\rho^2 |\nabla \theta|^2
{\rm d}x.
\end{equation}
Moreover, if we put
$$
\lambda=
\max\left\{\frac{9}{2\ve^2\min_{\overline G^\prime_\delta}{|\nabla \theta|^2}},
\, 2{d^\prime}^2\right\},
$$
under the additional assumption that $|w_t|\leq 2$ in
$G_\delta^\prime$, the following pointwise inequality
$2\ve^2\lambda|w_t-{\rm e}^{i\theta}|^2|\nabla \theta|^2\geq
\rho^2(|w_t|^2-1)^2$ holds in $G_\delta^\prime$ (see the proof of
Proposition \ref{testingmap}). Thus,
\begin{equation}
\label{dveotsenki}
L_\ve^{(d^\prime)}(w_t,G_\delta)\leq M^\prime_\lambda(w_t)+
\frac{1}{4\ve^2}\int_{G^{\prime\prime}_\delta}\rho^4(|w_t|^2-1)^2{\rm d}x.
\end{equation}

To demonstrate (\ref{dveotsenki1}) we first note that
\begin{equation}
\label{separMprime}
M^\prime_\lambda(w_t)=
%\frac{t^2\pi}{d^\prime} |b_{-1}(t)|^2\Phi^\prime_1(f_{-1})+
\frac{t^2\pi}{d^\prime}\sum_{k=-\infty}^\infty|b_k(t)|^2\Phi^\prime_k(f_k),
\end{equation}
where functionals $\Phi^\prime_k$ are defined
as functionals $\Phi_k$ in
(\ref{odnomernyi1}) with $d^\prime$ in place of $d$. The representation
(\ref{separMprime}) is analogous to that in (\ref{mmm}) while its justification
differs because of the fact that $\nabla h$ vanishes at least at some points of
the boundary of $G_\delta$ (and possibly somewhere in $G_\delta^{\prime\prime}$).
We rely on the following Lemma, which implies directly (\ref{separMprime}).

\begin{lem}
\label{separatedvariables}
Let $f,g\in C^1([1-\delta,1];\D{C})$ then for all integers $n,m$
$$
\int_{G_\delta} f(h){\rm e}^{in\theta}
\overline{g(h){\rm e}^{im\theta}}\, \rho^2|\nabla\theta|^2 {\rm d}x=
\begin{cases}
0,\quad \text{if}\quad n\not=m\\
\displaystyle
\frac{2\pi}{d^\prime}\int_{1-\delta}^1 f(s)\bar g(s){\rm d}s,\quad \text{if}\quad n=m.
\end{cases}
$$
\end{lem}
\noindent
{\bf Proof.}
By virtue of the pointwise equalities $\nabla h\cdot\nabla\theta=0$ and
$\div (\rho^2 \theta)=0$ in $G_\delta$, for any regular
values $\alpha,\beta$ of $h$ such that
$1-\delta\leq\alpha<\beta\leq 1$ and any integer $n\not=m$,
\begin{multline}
\label{o1}
\int_{\alpha<h<\beta}
f{\rm e}^{in\theta}\overline{g{\rm e}^{im\theta}}\, \rho^2|\nabla\theta|^2 {\rm d}x
=\frac{-i}{n-m}\int_{\alpha<h<\beta}\nabla \theta\cdot \nabla {\rm e}^{i(n-m)\theta} f\bar{g}\rho^2{\rm d}x\\
=\frac{i}{n-m}\int_{\alpha<h<\beta}
\div(\rho^2\nabla \theta)f\bar{g}{\rm e}^{i(n-m)\theta} {\rm d}x\\
+\frac{i}{n-m}\int_{\alpha<h<\beta}
\nabla\theta\cdot\nabla h(f^\prime\bar{g}+f\bar{g}^\prime){\rm e}^{i(n-m)\theta}
\rho^2{\rm d}x=0,
\end{multline}
where all integrals are understood over subsets of $G_\delta$.
If $n=m$ we set $F(h)=\int_{\alpha}^r f(s)\bar g(s)ds$ then,
since $|\nabla h|=d^\prime\rho^2|\nabla\theta|$
\begin{multline}
\label{o2}
\int_{\alpha<h<\beta} f(h)\bar{g}(h) \rho^2|\nabla\theta|^2{\rm d}x
=\frac{1}{{d^\prime}^2}\int_{\alpha<h<\beta}\nabla (F(h))\cdot \nabla h
\frac{{\rm d}x}{\rho^2}
=\frac{1}{{d^\prime}^2}F(\beta)\int_{h=\beta}
\frac{\partial h}{\partial \nu}\frac{{\rm d} s}{\rho^2}\\
-
\frac{1}{{d^\prime}^2}
\int_{\alpha<h<\beta}\div(\frac{1}{\rho^2}\nabla h)F(h) {\rm d}x=
\frac{2\pi}{d^\prime}\int_{\alpha}^\beta f(s)\bar{g}(s) {\rm d}s,
\end{multline}
here we have also used the fact that $\div(\frac{1}{\rho^2}\nabla h)=0$ in $G_\delta$.
The statement of the Lemma
is then obtained by passing to the
limits $\alpha\to 1-\delta$ and $\beta\to 1$ in (\ref{o1}, \ref{o2}).
\hfill $\square$

By using (\ref{asymtcoef}) in (\ref{separMprime}) we compute
\begin{multline}
\label{finalotsenka1}
M^\prime_\lambda(w_t)
=\frac{t^2\pi}{d^\prime}\sum_{k=-\infty}^\infty|c_k(t)|^2\Phi^\prime_k(f_k)+O(t^2)
=\pi((1-t)^2-1)^2\sum_{k=0}^\infty k(1-t)^{2k}\\
+2\pi t^2(\lambda-{d^\prime}^2)\sum_{k=1}^\infty \frac{(1-t)^{2k}}{k}+O(t^2)
=\pi-2\pi t+o(t).
\end{multline}
Explicit, but tedious computations (left to the reader), show also that
\begin{equation}
\label{finalotsenka2}
|w_t-{\rm e}^{id^\prime\theta}|\leq Ct\
\text{when}\
\theta\not \in (-\alpha,\alpha)+2\pi\D{Z},
\end{equation}
for any $0<\alpha<2\pi$, where $C$ is independent of $t$. From (\ref{finalotsenka2}), we
see the second term in (\ref{dveotsenki}) is of order $O(t^2)$ as $t\to 0$. Combined with (\ref{finalotsenka1}) this proves  (\ref{dveotsenki1}).

{\bf Final step.} The bound (\ref{dveotsenki}) is shown assuming that
$|w_t|<2$ in $G^{\prime\prime}_\delta$. Note, that this can always
be achieved by replacing $w_t$ by $\tilde w_t:=w_t \min\{1,2/|w_t|\}$, and this
change increases neither the first term in (\ref{dveotsenki}) nor the second one.
Thus in order to complete the proof of the lemma, we need
to show only that the map $v$ defined by (\ref{test}) satisfies
$d-1/2\leq{\rm abdeg}(v)\leq d+1/2$ when $t$ is chosen sufficiently
small. Indeed, due to (\ref{finalotsenka2}) $w_t$ weakly $H^1$-converges to
${\rm e}^{id^\prime\theta}$. Therefore the norm $\|u-v\|_{L^2(A)}$ tends
to $0$ when $t\to 0$. Then, according to Lemma \ref{cont}, for small $t$
${\rm abdeg}(v)$ is close to ${\rm abdeg}(u)$, while $d-1/2<{\rm abdeg}(u)<d+1/2$
and we are done.\hfill $\square$
%\newpage

\section*{Appendix A}

Here is a simple example of M$\rm\ddot{o}$bius (Blashke) test map
that illustrates an important property of nearboundary vortices.
Namely, each such vortex appears in pair with a "ghost"
antivortex (vortex with opposite sign) coming from outside the
domain. For simplicity consider the minimization problem of (\ref{P1})
with $A=B_1$, boundary condition $|u|=1$ on $\partial B_1$ and
imposed degree one on the boundary $\partial B_1$. It is shown in
\cite{BM1} that the infimum in this problem is not attained and
the behavior of the minimizing sequence is described by
%The behavior of a nearboundary
%vortex as it approaches the boundary is qualitatively described by
%the following well known M$\rm\ddot{o}$bius (Blashke)  map
\begin{equation*}
v_\ve(z)=\frac{\bar\zeta_\ve}{|\zeta_\ve|}
\frac{z-\zeta_\ve}{\bar\zeta_\ve z-1}=\frac{1}{|\zeta_\ve|}
\frac{z-\zeta_\ve}{z-(\zeta_\ve+o(1))},
\end{equation*}
where $\zeta_\ve\in B_1$, $|\zeta_\ve|\to 1$ as $\ve\to 0$.
%In
%view of it is convenient to chose
We assume that ${\rm
dist}(\zeta_\ve,\partial B_1)=o(\ve)$. Clearly this map has single
zero of degree one at $z=\zeta_\ve$. Introduce
$\zeta^*_\ve=1/\bar\zeta_\ve$, and write this map as
%\begin{multline}
$$
v_\ve(z)=
\frac{1}{|\zeta_\ve|}\frac{|z-\zeta_\ve|}{|z-\zeta^*_\ve|}
\frac{z-\zeta_\ve}{|z-\zeta_\ve|}\left(\frac{z-\zeta^*_\ve}{|z-\zeta^*_\ve|}\right)^{-1}
=\omega_{\ve}(z)\frac{z-\zeta_\ve}{|z-\zeta_\ve|}
\frac{\overline{z-\zeta^*_\ve}}{|z-\zeta^*_\ve|}, \eqno{\rm(A.1)}
$$
%\label{ghost}
%\end{multline}
where $\omega_{\ve}(z)$ is a real-valued function. Then, the last
factor in (A.1)
%(\ref{ghost})
corresponds to a ``ghost" vortex with the center at $\zeta^*_\ve\not\in B_1$.
The complex conjugation means this vortex has degree $-1$
which is why it is called antivortex.
Finally, $\omega_{\ve}(z)=1+o(1)$ outside the disk of radius
$\ve$ with center $z=\zeta_\ve$. These observations show
superposition of vortex and antivortex has almost no energy away
from the core, whereas, the contribution in the energy from outside the core
of an inner vortex is known to
grow as $\log(1/\ve)$.

\section*{Appendix B}

Let us show that, for any integers $p$, $q$, $d$, the set of
admissible testing maps in problem (\ref{constrainedminimization})
is not empty and
$$
m_\ve(p,q,d)\leq I_0(d,A)+\pi|q-d|+\pi|p-d|. \eqno{\rm(B.1)}
$$
To this end chose two sequences $(x_1^{(k)}), (x_2^{(k)})\subset A$
such that $x_1^{(k)}\to\partial\Omega$,
$x_2^{(k)}\to\partial\omega$. Consider solutions $\rho_1^{(k)}$,
$\rho_2^{(k)}$ of the boundary value problems
$$
\begin{cases}
\Delta \rho_1^{(k)}=2\pi|d-q|\delta(x-x_1^{(k)}) \ \text{in}\ A\\
\rho_1^{(k)}=0\quad \text{on}\ \partial \Omega,\\ \rho_1^{(k)}={\rm
Const} \quad \text{on}\
\partial\omega\\
\displaystyle\int_{\partial \Omega}
\frac{\partial\rho_1^{(k)}}{\partial \nu}{\rm d}\sigma=2\pi|d-q|,
\end{cases}
\ \text{\rm and}\
\begin{cases}
\Delta \rho_2^{(k)}=2\pi|d-p|\delta(x-x_2^{(k)}) \ \text{in}\ A\\
\rho_2^{(k)}=0\quad \text{on}\ \partial \omega,\\
\rho_2^{(k)}={\rm Const} \quad \text{on}\
\partial\Omega\\
\displaystyle\int_{\partial \omega}
\frac{\partial\rho_2^{(k)}}{\partial \nu}{\rm d}\sigma=2\pi|d-p|.
\end{cases}
$$
One can show that, for any neighborhoods $U$ and $V$ of
$\partial\Omega$ and $\partial\omega$ in $A$,
$$
\rho_1^{(k)}\to 0\ \text{\rm in}\ C^1(A\setminus U)\quad \text{\rm
and}\quad \rho_2^{(k)}\to 0\ \text{\rm in}\ C^1(A\setminus V).
\eqno{\rm(B.2)}
$$
There exist harmonic conjugates
$$
\begin{array}{l}
\psi_1^{(k)}:A\setminus \{x_1^{(k)}\}\to \D{R}\setminus
2\pi|d-q|\D{Z},\\
\psi_2^{(k)}:A\setminus \{x_2^{(k)}\}\to \D{R}\setminus
2\pi|d-p|\D{Z},
\end{array}
$$
(i.e. $\nabla\psi_j^{(k)}=\nabla^\bot\rho_j^{(k)}$, $j=1,2$) so that
$e^{\rho_j^{(k)}+i\psi_j^{(k)}}$ are harmonic functions in $A$. By
using the pointwise equality $\frac{1}{2}|\nabla u|^2=
\partial_{x_1} u\times \partial_{x_2} u+\frac{1}{4}|\partial_{\bar z} u|^2$
and integrating by parts we get
$$
\frac{1}{2}\int_A e^{2\rho_1^{(k)}}(|\nabla \rho_1^{(k)}|^2+|\nabla
\psi_1^{(k)}|^2)\,{\rm d}x=\frac{1}{2}\int_A |\nabla
e^{\rho_1^{(k)}+i\psi_1^{(k)}}|^2\,{\rm d}x=\pi|d-q|.
\eqno{\rm(B.3)}
$$
Similarly we have,
$$
\frac{1}{2}\int_A e^{2\rho_2^{(k)}}(|\nabla \rho_2^{(k)}|^2+|\nabla
\psi_2^{(k)}|^2)\,{\rm d}x=\pi|d-p|. \eqno{\rm(B.4)}
$$

Then we introduce
$$
u^{(k)}(x)=u(x) e^{\rho_1^{(k)}(x)+i\,{\rm sgn}(q-d)\psi_1^{(k)}(x)}
e^{ \rho_2^{(k)}(x) +i\,{\rm sgn}(p-d)\psi_2^{(k)}(x)}
e^{-\rho_1^{(k)}(\partial\omega)+
V(x)(\rho_1^{(k)}(\partial\omega)-\rho_2^{(k)}(\partial\Omega))},
$$
where $u$ is a minimizer of problem (\ref{S^1minimization}) and $V$
is the solution of (\ref{capacitysolution}). From (B.2-B.4) we
derive that, for $k$ sufficiently large, $u^{(k)}$ is an admissible
testing map for problem (\ref{constrainedminimization}), and that
$E_\ve(u^{(k)})\to I_0(d,A)+\pi|q-d|+\pi|p-d|$ when $k\to\infty$.

\end{document}